\documentstyle[11pt,bbold]{article}

{
\newenvironment{beweis}{{\it Proof.}\ }{\ $\ \ \ \Diamond$ \\\ }

\renewcommand{\labelenumi}{{$\left(\roman{enumi}\right)$}}

\newcounter{nsatz}[section]
\newcounter{nlemma}[section]
\newcounter{ndef}[section]
\newcounter{nhyp}[section]
\newcounter{nconjecture}[section]
\newcounter{ncor}[section]
\newcounter{nrem}[section]
\newcounter{nexample}[section]
\newcounter{nprop}[section]

\newenvironment{nsatz}{\refstepcounter{nsatz}{\bf \arabic{section}.\arabic{nsatz}}\ 
               {\sc\bf Theorem.\ }\it}{\\\\ \rm}

\newenvironment{nlemma}{\setcounter{nlemma}{\value{nsatz}}
               \refstepcounter{nlemma}
               \setcounter{nsatz}{\value{nlemma}}
               {\bf \arabic{section}.\arabic{nsatz}}\ 
               {\sc\bf Lemma.\ }\it}{\\\\ \rm}

\newenvironment{nconjecture}{\setcounter{nconjecture}{\value{nsatz}}\refstepcounter{nconjecture}
               \setcounter{nsatz}{\value{nconjecture}}
               {\bf \arabic{section}.\arabic{nsatz}}\ 
               {\sc\bf Conjecture.\ }\it}{\\\\ \rm}

\newenvironment{ncor}{\setcounter{ncor}{\value{nsatz}}
               \refstepcounter{ncor}
               \setcounter{nsatz}{\value{ncor}}
               {\bf \arabic{section}.\arabic{nsatz}}\ 
               {\sc\bf Corollary.\ }\it}{\\\\ \rm}

\newenvironment{nexample}{\setcounter{nexample}{\value{nsatz}}
               \refstepcounter{nexample}
               \setcounter{nsatz}{\value{nexample}}
               {\bf \arabic{section}.\arabic{nsatz}}\ 
               {\sc\bf Example.\ }}{\\\\ \rm}

\frenchspacing
\begin{document}
\setlength{\headheight}{-2.9cm}
\newcommand{\MengeZ}{{\mbox{{\sf Z\hspace{-0.4em}Z}}}}
\newcommand{\z}{{\mbox{{\sf Z\hspace{-0.4em}Z}}}}
\newcommand{\bigudot}{{ }^{\;\;\;{\bullet}}\!\!\!\!\!\!\!\;\!\bigcup}
{\renewcommand{\section}[1]{\refstepcounter{section}{\bf \thesection.}\ {\bf {\
#1}}}
\renewcommand{\r}{{\mbox{\rm I$\!$R}}}
\renewcommand{\labelenumi}{(\alph{enumi})}
\renewcommand{\labelenumii}{(\arabic{enumii})}
\renewcommand{\labelenumiii}{(\alph{enumiii})}
\mbox{\vspace{4cm}}
\vspace{4cm}
\begin{center}
{\bf \Large Fixed conjugacy classes of normal subgroups and the $k(GV)$--problem \\}
\vspace{3cm}
by\\
\vspace{11pt}
Thomas Michael Keller\\
Department of Mathematics\\
Texas State University\\  
601 University Drive\\
San Marcos, TX 78666\\
USA\\
e--mail: keller@txstate.edu\\
\vspace{1cm}
\end{center}
\thispagestyle{empty}
\renewcommand{\thefootnote}{\fnsymbol{footnote}}

\newcommand{\og}{\ \raisebox{-.3em}{$\stackrel{\scriptstyle \geq}
{\scriptstyle \sim}$} \ }
\newcommand{\ob}{o.B.d.A.\ }
\newcommand{\End}{\mbox{End}}
\newcommand{\GL}{\mbox{GL}}
\newcommand{\ug}{\ \raisebox{-.3em}{$\stackrel{\scriptstyle \leq}
{\scriptstyle \sim}$} \ }
\newcommand{\n}{{\mbox{\rm I$\!$N}}}       
\newcommand{\red}{$\sigma$--reduced}
\newcommand{\theorem}[1]{$   #1 \enspace {\mbox{\sc Theorem. }} $}
\newcommand{\rem}[1]{$  #1 \enspace {\mbox{\sc Remark. }} $}
\newcommand{\corollary}[1]{$  #1 \enspace {\mbox{\sc Corollary. }} $}
\newcommand{\defi}[1]{$  #1 \enspace {\mbox{\sc Definition. }} $}
\newcommand{\exa}[1]{$  #1 \enspace {\mbox{\sc Example. }} $}
\newcommand{\exas}[1]{$  #1 \enspace {\mbox{\sc Example. }} $}
\newcommand{\cor}[1]{$  #1 \enspace {\mbox{\sc Corollary. }} $}
\newcommand{\hyp}[1]{$  #1 \enspace {\mbox{\sc Hypothesis. }} $}
\newcommand{\lemma}[1]{$  #1 \enspace {\mbox{\sc Lemma. }} $}
\newcommand{\pr}{{\it  Proof. }  }
\newcommand{\notat}{{\mbox{\sc  Notation. }}  }
\newcommand{\tn}{\ \vert \hspace{-.73em} -}
\newcommand{\pp}{\mbox{\cal P}}
\newcommand{\X}{{\mbox{$\setminus$}\mbox{$\!\!\!/$}}}
\newcommand{\ind}{\not\ \!\!\!\!|}
\newcommand{\cd}{\mbox{\rm cd}}
\newcommand{\cl}{\mbox{\rm cl}}
\newcommand{\scl}{\mbox{\scriptsize\rm cl}}
\newcommand{\sirr}{\mbox{\scriptsize\rm Irr}}
\newcommand{\Irr}{\mbox{\rm Irr}}
\newcommand{\dl}{\mbox{\rm dl}}
\newcommand{\fdl}{\mbox{\rm\footnotesize dl}}
\newcommand{\tdl}{\mbox{\rm\tiny dl}}
\newcommand{\stab}{\mbox{\rm stab}}
\newcommand{\Gal}{\mbox{\rm Gal}}
\newcommand{\cp}{\raisebox{2pt}{\tiny{\ \sf Y}}}
\newcommand{\core}{\mbox{\rm core}}
\renewcommand{\char}{\mbox{\rm char}}
\newcommand{\Syl}{\mbox{\rm Syl}}
\setlength{\parindent}{0pt}
\newpage

\begin{center}
\parbox{12.5cm}{{\small
{\sc Abstract.}   
We establish several new bounds for the number of conjugacy classes of a finite 
group, all of which involve the maximal number  $c$ of conjugacy classes of a normal
subgroup fixed by some element of a suitable subset of the group. To apply these
formulas effectively, the parameter $c$, which in general is hard to control,
is studied in some important situations.\\
These
results are then used to provide a new, shorter proof of the most difficult case of the
well--known $k(GV)$--problem, which occurs for $p=5$ and $V$ 
induced from the natural module of a 
5--complement of $\GL (2,5)$. We also show how, for large $p$, the new 
results reduce the $k(GV)$--problem
to the primitive case, thereby improving previous work on this. Furthermore, 
we discuss how they can be used in tackling the 
imprimitive case of the as of yet unsolved noncoprime $k(GV)$--problem.\\
}}
\end{center}
\normalsize

\section{Introduction and notation}\label{section1}\\

Bounding the number of conjugacy classes of a finite group is a fundamental issue
in finite group theory, as is evidenced by the large body of literature on the subject
(for general results see e. g. \cite{gallagher}, \cite{kovacs-robinson}, 
\cite{liebeck-pyber}, for asymptotic results on classical groups see \cite{fulman-guralnick},
for the $k(GV)$--problem see e. g. \cite{robinson-thompson}, \cite{gmrs}). This paper
is another contribution to the subject, providing some general bounds involving a parameter
that, as far as we can tell, has hardly been used up to this point, 
but which will prove quite useful. \\
This new parameter
is $|C_{\scl (N)}(g)|$. Here $G$ is a group, $g\in G$, $N\unlhd G$, $\cl(N)$ is the set of conjugacy classes
of $N$, and $C_{\scl(N)}(g)$ is the set of classes of $N$ which are fixed (as a set) by 
$g$ under conjugation.  Note that by Brauer's permutation lemma
we have 
\[|C_{\scl (N)}(g)|\ =\ |C_{\sirr (N)}(g)|,\] 
where the latter is the number of 
irreducible complex characters of $N$ fixed by $g$.\\
This parameter shows up in a number of bounds for the number $k(G)$ of conjugacy classes of
$G$, such as the following:\\

{\bf Lemma A. }{\it
Let $G$ be a finite group and $N\unlhd G$. Let $g_i\in G$ ($i=1,\ldots,k(G/N)$) such that
the $g_iN$ are representatives of the conjugacy classes of $G/N$. Then
\[k(G)\ \leq\ \sum_{i=1}^{k(G/N)}|C_{\scl (N)}(g_i)|.\]
}
(See \ref{lem1.1} below.)\\

{\bf Lemma B. }{\it
Let $G$ be a finite group and $H\leq G$. Let $N$ be the core of $H$ in $G$. Then
\[k(G)\leq k(H)+k(G/N)\max\{|C_{\scl(N)}(g)|\ |\ g\in G-\bigcup_{x\in G}H^x\}.\]
}
(See \ref{lemb1} below.)\\

If one wants to use these and similar results effectively, one must somehow control
$|C_{\scl(N)}(g)|$, which seems to be very difficult in general, and we are not aware
of any result on this in the literature with the exception of our own first
encounter with it in \cite{kellerkgv1}, where some very technical result on it
was proved in \cite[Lemma 4.7(b)]{kellerkgv1}.). 
In Sections \ref{sectiona}, \ref{sectionc}, and \ref{sectiond},
therefore, we will prove some bounds on $|C_{\scl(N)}(g)|$ in some key special cases.\\

These techniques, while technical at times, turn out to be quite powerful. We will
demonstrate this in Section \ref{sectione}, where we will give a short proof of the most
difficult case of the $k(GV)$--problem that has only been solved recently by the combined
efforts of Gluck, Magaard, Riese, and Schmid in \cite{gmrs}. Recall that the $k(GV)$--problem claims that
$k(GV)\leq |V|$ whenever $V$ is a finite faithful $G$--module of characteristic $p$ not
dividing $|G|$. This problem, which is equivalent to Brauer's well--known $k(B)$--problem
for $p$--solvable groups, has kept mathematicians busy for the past 20 years, and the final
step in its solution was a special case for $p=5$ treated in \cite{gmrs} that had
escaped all former attacks. So we will give a new proof of this case. More precisely, we
will prove (see \ref{theoreme3} below):\\

{\bf Theorem C. }{\it
Let $G$ be a finite $5'$--group and $V$ be a fiathful GF$(5)$--module such that $V$ is induced
from a $G_1$--module $W$, where $G_1$ is a suitable subgroup of $G$, $|W|=25$ and
$G_1/C_{G_1}(W)\not=1$ is isomorphic to a subnormal subgroup of $L$, where $L$ is a 5--complement
in $\GL (2,5)$. Suppose that whenever $U\leq G$ and $X\leq V$ is a $U$--module with $|UX|<|GV|$, 
then $k(UX)\leq|X|$.
Then
\[k(GV)\leq|V|.\]
}

In \ref{sectiond} we will also use our techniques to directly 
reduce the $k(GV)$--problem to the case of $V$
being primitive as $G$--module whenever $p>2^{47}$ (see \ref{lemd4}). This improves
and shortens the corresponding reduction in our previous proof of the $k(GV)$--problem for
large primes in \cite[Theorem 4.1]{kellerkgv1}. Finally, in Section \ref{sectionf} we turn
to the more recent non--coprime $k(GV)$--problem (see \ref{conf1}), research on which
is still in its beginnings. While some work by Guralnick and Tiep \cite{guralnick-tiep}
on primitive groups is underway, nothing is known on how to deal with the imprimitive
case. We provide a few first steps in this direction that might be useful in an inductive
argument (see \ref{corf3} and \ref{theorem7.4} and the remarks following each of them). For instance, we will prove
the following.\\

{\bf Theorem D. }{\it
Let $G$ be a finite group and $V$ be a finite $G$--module. Suppose that $N\unlhd G$ and
$V_N=V_1\oplus\ldots\oplus V_n$ for some $n\geq 5$ such that $G/N$ primitively and faithfully
permutes the $V_i$. Moreover suppose that for some constant $C>0$ we have
\[k(NV)\leq C|V|\log_2|V|\]
and $k(UV_1)\leq C|V_1|\log_2|V_1|$ for every $U\leq N/C_N(V_1)$, and
\[|N/C_N(V_1)|\leq\frac{1}{50}C^{\frac{14}{3n}-\frac{8}{3}}|V_1|(\log_2|V|)^{\frac{14}{3n}-\frac{8}{3}}.\]
If $F^*(G/N)$ (the generalized Fitting subgroup of $G/N$) is not a product of alternating
groups, then
\[k(GV)\leq C|V|\log_2|V|.\]
}

Our notation is as in \cite{kellerkgv1} and \cite{kellerkgv2}. In particular, if $G$
acts on a set $\Omega$, we denote by $n(G,\Omega)$ the number of orbits of $G$ on
$\Omega$ and by $C_\Omega(g)$ the set of fixed points of $g$ on $\Omega$. We will freely
use the elementary formulas for $k(GV)$ as discussed in \cite{kellerkgv1} as well as 
the well--known fact that if $N\unlhd G$, then $k(G)\leq (G/N)k(N)$. We will also
use the latest improvement on upper bounds for the number of conjugacy classes of permutation
groups. This is due to A. Mar\'{o}ti \cite[Theorem 1.1]{maroti} and states that for every $U\leq S_n$ 
with $n\ne 2$ we have $k(U)\leq 3^{(n-1)/2}$.\\

\section{On conjugacy classes fixed by an automorphism}\label{sectiona}

In this section we study the action of group elements on the conjugacy classes of some normal
subgroup of the group.\\
Bits and pieces of what is to follow have already been foreshadowed in 
\cite{kellerkgv1} and \cite{kellerkgv2}, but our treatment here is
self--contained.\\

We start with a general lemma.\\

\begin{nlemma}\label{lem1.1}
Let $G$ be a finite group and $N\unlhd G$. Let $g_iN$ ($i=1,\ldots,k(G/N)$)
be representatives of the conjugacy classes of $G/N$, and write $\Omega_i$ for the set
of $N$--orbits on $g_iN$. Then
\[k(G)\ =\ \sum_{i=1}^{k(G/N)}n(C_{G/N}(g_i),\Omega_i)\ \leq\ \sum_{i=1}^{k(G/N)}|C_{\scl (N)}(g_iN)|\]
\end{nlemma}
\begin{beweis}
Let $\Omega=\{g^N|g\in G\}$ be the set of $N$--orbits of $G$, and for 
$g\in G$ put $\Omega_{gN}=\{\omega\in\Omega|\omega\subseteq gN\}$.
Hence $\Omega_i=\Omega_{g_iN}$ for $i=1,\ldots,k(G/N)$.
For each $\omega\in\Omega$ let $g_\omega\in G$, so that $g_\omega^N=\omega$. 
It is easy to check that $C_{G/N}(\omega)\leq C_{G/N}(g_\omega N)$ and that for 
$\omega\in\Omega_{gN}$ we have $g_\omega N=gN$.
With this we conclude that
\begin{eqnarray*}
k(G)&=&n(G/N,\Omega)=\frac{1}{|G/N|}\sum_{\omega\in\Omega}|C_{G/N}(\omega)|\\
    &=&\frac{1}{|G/N|}\sum_{\omega\in\Omega}|C_{C_{G/N}(\omega)}(g_\omega N)|\\
    &=&\frac{1}{|G/N|}\sum_{gN\in G/N}\sum_{\omega\in\Omega_{gN}}|C_{G/N}(\omega)\cap
       C_{G/N}(g_\omega N)|\\
    &=&\frac{1}{|G/N|}\sum_{gN\in G/N}|C_{G/N}(gN)|\frac{1}{|C_{G/N}(gN)|}
       \sum_{\omega\in\Omega_{gN}}|C_{C_{G/N}(gN)}(\omega)|\\
    &=&\frac{1}{|G/N|}\sum_{gN\in G/N}|C_{G/N}(gN)|n(C_{G/N}(gN),\Omega_{gN})\\
    &=&\sum_{i=1}^{k(G/N)}n(C_{G/N}(g_iN),\Omega_i)\qquad (*)
\end{eqnarray*}
which is the first part of the lemma.
(Notice that $(*)$ was already proved in \cite[Lemma 1.6]{kellerkgv1} with a much
longer proof.)\\
Next fix $g\in G$. We claim that $(**)$ $|\Omega_{gN}|=|C_{\Omega_N}(g)|$,
that is, the number of $N$--orbits on $gN$ is the same as the number of
conjugacy classes of $N$ fixed by $g$. To see this, let $x\in N$. Then 
$C_{gN}(x)$ is nonempty if and only if there is an $n\in N$ with $x^{-1}gnx=gn$
which means that $x^{n^{-1}}=x^g$. This is equivalent to $x^N\in C_{\Omega_N}(g)$.
Moreover then obviously $C_{gN}(x)=gnC_N(x)$. So altogether it follows that 
$|C_{gN}(x)|=|C_N(x)|$ whenever $C_{gN}(x)\not=\emptyset$. Keeping all this in mind
and putting $X=\{y\in N|y^N\in C_{\Omega_N}(g)\}$, we see that $y\in X$ if and only if 
$C_{gN}(y)\not=\emptyset$ and thus
\[|C_{\Omega_N}(g)|\ =\ n(N,X)\ =\ \frac{1}{|N|}\sum_{y\in X}|C_N(y)|\ =\ 
\frac{1}{|N|}\sum_{y\in X}|C_{gN}(y)|\ =\ \frac{1}{|N|}\sum_{y\in N}|C_{gN}(y)|
\ =\ |\Omega_{gN}|,\]
where the last equation follows with the Cauchy--Frobenius orbit counting formula.
So $(**)$ is proved.\\
Now as clearly $\Omega_N=\cl(N)$ and $C_{\Omega_N}(g)=C_{\Omega_N}(gN)$,
by $(**)$ we conclude that $n(C_{G/N}(g_iN),\Omega_i)\leq |\Omega_i|=
|C_{\scl(N)}(g_iN)|$, and so by $(*)$ the assertion of the lemma follows.
\end{beweis}

\begin{nlemma}\label{lem1.2}
Let $G$ be a finite group. Suppose that $M\unlhd G$ and that $M=\X_{i=1}^l M_i$, where
the $M_i$ ($i=1,\ldots,l$) are normal subgroups of $G$. We write elements $(a_1,\ldots,a_l)\in M$
simply as $a_1\ldots a_l$ (for $a_i\in M_i$).
Moreover suppose that
$G/M=\langle gM\rangle$ is cyclic of order $m$. Let $N\leq M$ with $N^g=N$, and put
$L_i=M_i\times\ldots\times M_l$ for $i=1,\ldots,l+1$ (so $L_{l+1}=1$).\\
(a) Let $x=x_1\cdot\ldots\cdot x_l\in N$, where $x_i\in M_i$ ($i=1,\ldots,l$),
and
put $C_i=\bigcap\limits_{j=1}^{i-1}C_N(x_j)$ for $i=2,\ldots,l$, and put $C_1=N$.\\
Then the following are equivalent:
\begin{itemize}
\item[(i)]$x^g\in x^N$
\item[(ii)]For $i=1,\ldots,l$ there exist $z_i\in C_i$ and
           $gz_1\cdot\ldots\cdot z_i\in C_G(x_i)$.
\item[(iii)]Put $K_i=\{y\in M_i\ |\ x_1\cdot\ldots\cdot x_{i-1}yL_{i+1}\in NL_{i+1}/L_{i+1}\}$
            for $i=1,\ldots,l$. Note that $C_i$ acts on $K_i$ (by conjugation). For 
            $i=1,\ldots,l-1$ there exist $z_i\in C_i$ such that 
            $x_i^{gz_1\cdot\ldots\cdot z_{i-1}}\in K_i$ and 
            $x_i^{gz_1\cdot\ldots\cdot z_{i-1}}$ and $x_i$ lie in the same orbit of 
            $C_i$ on $K_i$. (The $z_i$ here are actually the same as in (ii).)
\end{itemize}

(b) Let $k_i=\max\{|C_{\scl(U)}(h)|\ |\ h\in G-M\mbox{ and }U\leq M_i\mbox{ with }
U^h=U\mbox{ and }h^m\in M_1\times\ldots\times M_{i-1}\times U\times M_{i+1}\times\ldots\times M_l\}$ 
for $i=1,\ldots,l$. Then
\[|C_{\scl(N)}(g)|\ \leq\ \prod_{i=1}^lk_i.\]
\end{nlemma}
\begin{beweis}
(a) (i) $\Rightarrow$ (ii): Put $x_0=z_0=1$ and $C_0=G$. We show by induction
on $i$ that there are $z_i\in C_i$ ($i=0,\ldots,l$) such that 
$gz_0z_1\cdot\ldots\cdot z_i\in C_G(x_i)$. For $i=0$ this is trivial. So let
$1\leq i\leq l$ and suppose that we already have $z_0,\ldots,z_{i-1}$.
Then $x^{gz_0\cdot\ldots\cdot z_{i-1}}=\prod\limits_{j=1}^lx_j^{gz_0\cdot\ldots\cdot z_{i-1}}=
x_1\cdot\ldots\cdot x_{i-1}\prod\limits_{j=i}^lx_j^{gz_0\cdot\ldots\cdot z_{i-1}}$.
As $x^g\in x^N$, we know that there is a $z_i\in N$ such that 
$x=x^{(gz_0\cdot\ldots\cdot z_{i-1})z_i}=x_1^{z_i}\cdot\ldots\cdot x_{i-1}^{z_i}
\prod\limits_{j=i}^lx_j^{gz_0\cdot\ldots\cdot z_i}$. This clearly forces 
$z_i\in C_i$ and $gz_0\cdot\ldots\cdot z_i\in C_G(x_i)$, as wanted.\\

(ii) $\Rightarrow$ (iii): Let the $z_i\in C_i$ ($i=1,\ldots,l$) be as in
(ii) and fix $i\in\{1,\ldots,l-1\}$. Then $x^{gz_1\ldots z_{i-1}}L_{i+1}=
x_1\cdot\ldots\cdot x_{i-1}x_i^{gz_1\ldots z_{i-1}}L_{i+1}\in NL_{i+1}/L_{i+1}$
which shows that $x_i^{gz_1\ldots z_{i-1}}\in K_i$.
Moreover as $z_i\in C_i$ and $\left(x_i^{gz_1\ldots z_{i-1}}\right)^{z_i}=x_i$
we see that $x_i^{gz_1\ldots z_{i-1}}$ lie in the same orbit of $C_i$ on $K_i$.\\

(iii) $\Rightarrow$ (i): Suppose we have $z_i\in C_i$ ($i=1,\ldots,l-1$) as in
(iii). As $x_l^{gz_1\ldots z_{l-1}}$ and $x_l$ lie in the same orbit of $C_l$ on
$K_l$, there is a $z_l\in C_l$ such that $x_l^{gz_1\ldots z_l}=x_l$.
Thus clearly $x^{gz_1\ldots z_l}=x$, and so $x^g=x^{(z_1\ldots z_l)^{-1}}\in x^N$.\\

(b) We prove the statement by induction on $l$. If $l=1$, the statement is easily seen
to be true. Let $l\geq 1$. Observe that 
$M_1\times\ldots\times M_{l-1}\cong M_0:=M/M_l\unlhd G/M_l=:G_0$ and consider $N_0=NM_l/M_l$. Then 
by induction we have
\[|C_{\scl(N_0)}(gM_l)|\ \leq\ \prod_{i=1}^{l-1}k_{0i},\mbox{ where}\]
\begin{eqnarray*}
k_{0i}&=&\max\{|C_{\scl(U)}(hM_l)|\ |\ hM_l\in G_0-M_0\mbox{ and }U\leq M_iM_l/M_l
         \mbox{ with }U^{hM_l}=U\mbox{ and }\\
       &&(hM_l)^m\in M_1M_l/M_l\times\ldots\times M_{i-1}M_l/M_l\times U\times M_{i+1}M_l/M_l
          \times\ldots\times M_{l-1}M_l/M_l\},
\end{eqnarray*}
and as $M_iM_l/M_l\cong M_i$ (as $G$--sets) and $M_l$
centralizes $M_i$ for $i<l_i$ we see that $k_{0i}=k_i$ for $i=1,\ldots,l-1$.\\
Now if $x_i\in M_i$ ($i=1,\ldots,l$) such that for $x=x_1\ldots x_l$ we have
$x^g\in x^N$, then $xM_l=x_1\ldots x_{l-1}M_l$ satisfies $(xM_l)^{gM_l}\in (xM_l)^{N_0}$,
and so by the above $xM_l$ is in one of at most $\prod\limits_{i=1}^{l-1}k_i$ possible
conjugacy classes of $N_0$.\\
Next suppose that $xM_l\in N_0$ with $(xM_l)^{g}\in (xM_l)^{N_0}$ has been chosen, i. e.,
$x_i\in M_i$ ($i=1,\ldots,l-1$) are already fixed. Let $K_l,C_l$ be as in (a).
Now if $y_1,\ y_2\in K_l$, then
$x_1\ldots x_{l-1}y_1$ and $x_1\ldots x_{l-1}y_2$ obviously are in the same class
of $N$ if and only if $y_1$ and $y_2$ are in the same orbit of $C_l$ on $K_l$,
and thus there are exactly $n(C_l,K_l)$ different classes $y^N$ of $N$ such that
$(yM_l)^{N_0}=(xM_l)^{N_0}$. If $y^N$ is such a class and we choose the representative
$y\in N$ such that $y=x_1\ldots x_{l-1}x_l$ for a suitable $x_l\in K_l$, then by (a)
$y^N$ is fixed by $g$ only if $x_l$ and $x_l^h$ lie in the same orbit of $C_l$
on $K_l$ for some element $h\in C_G(x_1\ldots x_{l-1})\cap gN$, i. e., 
\[y^N\in C_{\scl(N)}(g)\mbox{ only if }(x_l^{C_l})^h=x_l^{C_l}\]
(note that $h$ clearly normalizes $C_l=C_N(x_1\ldots x_{l-1})$).
Now put 
\[L_l\ =\ \{z\in M_l\ |\ az\in C_l\mbox{ for some }a\in M_1\times\ldots\times M_{l-1}\}.\]
Observe that $L_l\leq M_l$ 
and that $K_l\subseteq L_l$, in particular $x_l\in L_l$.
Moreover, clearly $x_l^{C_l}=x_l^{L_l}$ and $L_l^h=L_l$ and $h^m\in M_1\times\ldots\times
M_{l-1}\times L_l$. Thus we see that
$y^N\in C_{\scl (N)}(g)$ only if $x_l^{L_l}\in C_{\scl (L_l)}(h)$, in which case
clearly $x_1\ldots x_{l-1}x_l^{L_l}\subseteq y^N$. Thus, 
since obviously $|C_{\scl (L_l)}(h)|\leq k_l$,
altogether we obtain
\[|C_{\scl(N)}(g)|\leq\left(\prod\limits_{i=1}^{l-1}k_i\right)\cdot k_l,\]
and so we are done.
\end{beweis}

\begin{nlemma}\label{lem2}
Let $G$ be a finite group. Suppose that $M\unlhd G$ and that $M=\X_{i=1}^pM_i$,
where $p$ is a prime and the $M_i$ are subgroups of $G$.
Moreover, suppose that $G/M=\langle gM\rangle$ is cyclic of order $p$, and that
$g$ permutes the $M_i$ transitively.\\
Let $N\leq M$ with $N^g=N$. Put $L=\{y\in M_1|ya\in N\mbox{ for some }a\in M_2\times
\ldots\times M_p\}$. (Clearly $L\leq M_1$.) Then
\[|C_{\scl(N)}(g)|\ \leq\ |L|.\]
\end{nlemma}
\begin{beweis}
Put $H=\langle N,g\rangle\leq G$. Then both $N$ and $H$ act on $N$ by conjugation, and
$k=|C_{\scl (N)}(g)|$ is exactly the number of common orbits of $H$ and $N$ on $N$. It is
an elementary fact (see \cite[Lemma 3.1]{guralnick-wan}) that
\[k\ =\ \frac{1}{|N|}\sum_{n\in N}|C_N(gn)|.\]
As $gn$ permutes the $M_i$ transitively (for each $n\in N$), it is clear that
$|C_N(gn)|\leq |L|$ for all $n\in N$. So the assertion follows.
\end{beweis}

For $p=2$, the above result turns out to be a little too weak for our purposes, and 
so later we will have to do some extra work to get around this.
The bound provided by \ref{lem2} is crude at times, and it is tempting to
believe that
\[|C_{\scl(N)}(g)|\leq\max\{k(U)|U\leq M_1\}\]
holds. This is not true, however, as the following example shows:\\

\begin{nexample}\label{ex0.3a}
Let $G=S_3\wr C_2$, where $C_2=\langle g\rangle$, and let 
$N=(S_3'\times S_3')\cdot\langle((12),(12))\rangle\leq S_3\times S_3$ (where $S_3'$ is the
commutator subgroup of $S_3$). So $|N|=18$, and $g$ normalizes $N$. It then can easily be
checked that $k(N)=6$ and $|C_{\scl(N)}(g)|=4$, so in particular
$|C_{\scl(N)}(g)|>\max\{k(U)|U\leq S_3\}=3$.\\
(Note that \ref{lem2} yields $|C_{\scl(N)}(g)|\leq 6$ here.)
\end{nexample}

In some sense, this seems to be a rare example, depending on the prime 2, as we see when we try
to use the above example to create a more general one: Let $q,p$ be primes with
$q|p-1$ and let $F$ be the Frobenius group of order $qp$. Put $G=F\wr C_2$, $C_2=\langle g\rangle$, 
and let, as above, $N$ be of order $p^2q$ such that $g$ normalizes $N$. Then one can check (by hand) that
if $q>2$, then $|C_{\scl(N)}(g)|=1+\frac{p-1}{q}+q-1=k(F)$, so here $|C_{\scl(N)}(g)|\leq
\max\{k(U)|U\leq F\}$, and if $q=2$, then 
\[|C_{\scl(N)}(g)|=1+2\frac{p-1}{q}+q-1=p+1,\]
so
\[|C_{\scl(N)}(g)|>\max\{k(U)|U\leq F\}=p.\]

\section{The general tools}\label{sectionb}\\

We now present our inductive arguments for proving results on $k(G)$.\\

\begin{nlemma}\label{lemb1}
Let $G$ be a finite group and $H\leq G$. Put $N=\bigcap\limits_{g\in G}H^g\unlhd G$. Then
\[k(G)\ \leq\ k(H)\ +\ k_0(G/N)\max\{|C_{\scl (N)}(g)|\ |\ g\in G-\bigcup\limits_{x\in G}H^x\},\]
where $k_0(G/N)$ is the number of conjugacy classes of $G/N$ that are contained in 
$G/N-\bigcup\limits_{x\in G}(H/N)^x$. (In particular, $k_0(G/N)\leq k(G/N)$.)
\end{nlemma}
\begin{beweis}
Note that by \ref{lem1.1} we have 
\[k(G)=\sum\limits_{i=1}^{k(G/N)}n(C_{G/N}(g_iN),\Omega_i),\] 
where the $g_iN$ are representatives of the
conjugacy classes of $G/N$, and $\Omega_i$ is the set of $N$--orbits on $g_iN$
(where $N$ acts by conjugation). Now let 
\[T\ =\ \{(gN)^{G/N}\ |\ g\in H\}\] 
be the set of conjugacy 
classes of $G/N$ that intersects nontrivially with $H/N$ and clearly we may assume that
\[T\ =\ \{(g_iN)^{G/N}\ |\ i=1,\ldots,|T|\}\] 
and that $g_i\in H$ for $i=1,\ldots,|T|$, and that
$g_1=1$. Then we have
\[k(G)\ =\ \sum_{i=1}^{|T|}n(C_{G/N}(g_iN),\Omega_i)\ +\ \sum_{i=|T|+1}^{k(G/N)}
n(C_{G/N}(g_iN),\Omega_i).\]
The second sum is bounded above by
\begin{eqnarray*}
    &&k_0(G/N)\max\{|\Omega_i|\ |\ i=|T|+1,\ldots,k(G/N)\}\\
    &&\leq k_0(G/N)\max\{|C_{\scl (N)}(g_iN)|\ |\ i=|T|+1,\ldots,k(G/N)\}\\
    &&\leq k_0(G/N)\max\{|C_{\scl (N)}(g)|\ |\ g\in G-\bigcup_{x\in G}H^x\}
\end{eqnarray*}
where the first inequality follows from the proof of \ref{lem1.1}. Thus it remains to show that
\[\sum\limits_{i=1}^{|T|}n(C_{G/N}(g_iN),\Omega_i)\leq k(H).\]
Now let $h_i\in H$, $i=1,\ldots,k(H/N)$,
be such that the $h_iN$ are representatives of the conjugacy classes of $H/N$ and let
$\Sigma_i$ be that set of $N$--orbits on $h_iN$ (with respect to conjugation). Clearly we may
assume that $h_i=g_i$ for $i=1,\ldots,|T|$. Then we conclude that 
\[\sum\limits_{i=1}^{|T|}n(C_{G/N}(g_iN),\Omega_i)\ =\ \sum_{i=1}^{|T|}n(C_{G/N}(h_iN),\Sigma_i)
\ \leq\ \sum_{i=1}^{k(H/N)}n(C_{H/N}(h_iN),\Sigma_i)\ =\ k(H),\]
where the last equality again follows from \ref{lem1.1}.
Hence the lemma is proved.
\end{beweis}

Note that \ref{lemb1} always yields $k(G)>k(H)$ which, for typical applications where
$k(H)$ is bounded by an inductive hypothesis, may give too weak a result unless additional
information is known on $k(H)$. We therefore also present another lemma that is more specialized, 
but possibly more suitable for inductive
arguments. This actually is a generalization of \cite[Lemmas 2.2 and 2.3]{kellerkgv2}.
The proof is quite similar to the proofs of those lemmas, but for the convenience of the reader
we outline the full argument.\\

\begin{nlemma}\label{lemb2}
Let $G$ be a finite group and $V$ be a finite $G$--module. Suppose that $N\unlhd G$ and
$V_N=V_1\oplus\ldots\oplus V_n$ for an $n\in\n$, where the $V_i$ are $N$--modules, and
assume that $G/N$ transitively and faithfully permutes the $V_i$. Put $H=N_G(V_1)$ and let
$f:\n\to\r$ be a function. Let $W\leq V$ be a $G$--submodule with $|W|\geq|V|^\delta$ for some
$0\leq\delta\leq 1$. Assume that there is a $G$--module $W'\leq V$ such that $V=W\oplus W'$.
Put
\begin{eqnarray*}
m_0&=&\max\{|C_{\scl(NV)}(g)|\ |\ g\in G, g\mbox{ has at most }\frac{n}{2}\mbox{ fixed points}\\  
    &&\mbox{in its permutation action on }\{V_1,\ldots,V_n\}\},
\end{eqnarray*}
and suppose that the following hold:
\begin{itemize}
\item[(i)] $k(HW)\leq f(|W|)$
\item[(ii)] $k(UN/N)\leq\frac{1}{\sqrt{n+1}}\left(\frac{f(|W|)}{m_0}\right)^\frac{1}{2}$
      for all $U\leq G$.
\end{itemize}
Then $k(GW)\leq f(|W|)$.
\end{nlemma}
\begin{beweis}
We may assume that $n\geq 2$. We consider the action of $G/N$ on $\Omega:=\Irr (NW)$. 
If $\omega\in\Omega$, we will write $\omega^G$ for the orbit of $\omega$ under $G$ and
$\omega\uparrow^G$ for the induced character. Let $P$ be Gallagher's goodness property with
respect to this action (see \cite[Example 3.4(b)]{kellerkgv1}. Then we have $k(HW)=\alpha_P(H/N,\Omega)$
and $k(GW)=\alpha_P(G/N,\Omega)$.
Now let 
\[R=\{gN\in G/N|gN\mbox{ normalizes at least }\frac{n}{2}\mbox{ of the }V_i\},\] 
so $R$ is a normal subset of $G/N$. Let 
\[T=\{\omega^{G/N}|\omega\in\Omega\mbox{ and }C_{G/N}(\omega)\not\subseteq R\},\] 
so $\omega^{G/N}$ means that there is a $gN\in G/N-R$ 
such that $gN$ fixes an element of $\omega^{G/N}$, i.e., $\omega^{G/N}\cap C_\Omega(gN)\not=\emptyset$.
Hence $\omega^{G/N}\cap C_\Omega(g^hN)\not=0$ for all $h\in G$. This shows that if $g_iN$,
$i=1,\ldots,t$ are representatives of the conjugacy classes of $G/N$ which are not in $R$, then
\[T\ \subseteq\ \bigcup\limits_{i=1}^t\{\omega^{G/N}\ |\ \omega\in\Omega\mbox{ and }\omega^{G/N}\cap
C_\Omega(g_iN)\not=\emptyset\}\]
and thus 
\[|T|\leq t\cdot\max\limits_{i=1,\ldots,t}|C_\Omega(g_iN)|\leq k(G/N)\max\limits_{gN\in(G/N)-R}
|C_\Omega(gN)|.\] 
Now if $gN\in G/N-R$, then $gN$ has at most $\frac{n}{2}$ fixed points on its 
permutation action on $\{V_1,\ldots,V_n\}$. Hence if we put $\Omega_1=\cl(NV)$, then
we have $|C_{\Omega_1}(g)|\leq
m_0$. Let $\Omega_0=\Irr(NV)$. Since by hypothesis there is a $G$--module $W'$ such that
$V=W\oplus W'$, clearly $NW\cong NV/W'$ and hence $\Omega\subseteq\Omega_0$. Now $G/N$ acts on 
both $\Omega_1$ and $\Omega_0$ by conjugation, and so Brauer's permutation lemma (see e.g.
\cite[Theorem 18.5(b)]{huppertcharacters}) yields 
\[|C_\Omega(gN)|\ \leq\ |C_{\Omega_0}(gN)|\ =\ |C_{\Omega_1}(gN)|\ \leq\ m_0 .\]
Hence we conclude that 
\[|T|\leq k(G/N)m_0.\]
Now consider $\omega$ with $\omega^{G/N}\not\in T$. Then $C_{G/N}(\omega)\subseteq R$,
so all elements of $C_G(\omega)N/N$ have at least $\frac{n}{2}$ fixed points
on $\{V_1,\ldots,V_n\}$. By \cite[Lemma 2.1]{kellerkgv2} there is an $i\in\{1,\ldots,n\}$ with
$C_G(\omega)\leq N_G(V_i)$, and so we may assume that $C_G(\omega)\leq H$.
As $H<G$, it follows that $|\omega^{G/N}|>|\omega^{H/N}|$, and so if
$\omega_1,\ldots,\omega_k\in\omega^{G/N}$ are representatives of the orbits
of $H/N$ on $\omega^{G/N}$ with $\omega_1=\omega$, then $k\geq 2$, and by the
Theorem in \cite{gallagher} and \cite[Exercise E17.2]{huppertcharacters} we see
that for $i=2,\ldots,k$ we have
\begin{eqnarray*}
k_P(C_{H/N}(\omega))&=&k_P(C_{G/N}(\omega))\ =\ k_P(C_{G/N}(\omega_i))\\
                    &=&|\{\psi\in\Irr (C_{GV}(\omega_i))\ |\ \psi\mbox{
                       is a constituent of the induced}\\
                       &&\mbox{ character }\omega_i
                       \uparrow^{
                       C_{GV}(\omega_i)}\}|\\
                    &\leq&|C_{GV}(\omega_i):C_{HV}(\omega_i)|\cdot
                          |\{\Theta\in\Irr (C_{HV}(\omega_i))\ |\ \Theta\mbox{
                       is a}\\
                       &&\mbox{ constituent of }\omega_i
                       \uparrow^{
                       C_{HV}(\omega_i)}\}|\\
                    &\leq&|G:H|\cdot k_P(C_{H/N}(\omega_i))\\
                    &=&n\cdot k_P(C_{H/N}(\omega_i)),      
\end{eqnarray*}
and hence we obtain that
\begin{eqnarray*}
\sum_{j=1}^kk_P(C_{H/N}(\omega_j))&\geq& k_P(C_{H/N}(\omega_1))+(k-1)\frac{1}{n}k_P(C_{H/N}(\omega))\\
                                  &\geq&\frac{n+k-1}{n}k_P(C_{H/N}(\omega))\\
                                  &\geq&\frac{n+1}{n}k_P(C_{H/N}(\omega))
\end{eqnarray*}
Since these considerations hold for any $\omega^{G/N}\not\in T$, we conclude 
that if $\omega_i\in\Omega$ ($i=1,\ldots,n(G/N,\Omega)$) are representatives
of the orbits of $G/N$ on $\Omega$ and the $\omega_{ij}$ ($j=1,\ldots,k_i$) are
representatives of the orbits of $H/N$ on $\omega_i^{G/N}$,
then we may assume that for all $i$ with $\omega_i^{G/N}\not\in T$ we have
$C_{G/N}(\omega_i)\leq H/N$, and then the above yields
\begin{eqnarray*}
\sum_{i\ \mbox{\scriptsize with }\omega_i^{G/N}\not\in T}k_P(C_{G/N}(\omega_i))&=&
\sum_{i\ \mbox{\scriptsize with }\omega_i^{G/N}\not\in T}k_P(C_{H/N}(\omega_i))\
\leq \frac{n}{n+1}\sum_{i\ \mbox{\scriptsize with }\omega_i^{G/N}\not\in T}
      \sum_{j=1}^{k_i}k_P(C_{H/N}(\omega_{ij}))\\
&\leq&\frac{n}{n+1}\alpha_P(H/N,\Omega)\
=\ \frac{n}{n+1}k(HW).     
\end{eqnarray*}
Hence altogether we obtain
\begin{eqnarray*}
k(GW)&=&\alpha_P(G/N,\Omega)\\
     &=&\sum_{i\ \mbox{\scriptsize with }\omega_i^{G/N}\in T}
         k_P(C_{G/N}(\omega_i))\ +\ \sum_{i\ \mbox{\scriptsize with }\omega_i^{G/N}\not\in T}k_P(C_{G/N}(\omega_i))\\
    &\leq&|T|\cdot\max_{i=1,\ldots,n(G/N,\Omega)}k_P(C_{G/N}(\omega_i))\ +\ \frac{n}{n+1}k(HW)\\
    &\leq&k(G/N)m_0\max_{U\leq G}k(UN/N)\ +\ \frac{n}{n+1}k(HW)\\
     &\leq&\left(\max_{U\leq G}k(UN/N)\right)^2\cdot m_0\ +\ \frac{n}{n+1}k(HW).
\end{eqnarray*}
Therefore by our Hypotheses (i) and (ii) we are done.
\end{beweis}

The final lemma in this section will be useful in certain noncoprime situations.\\

\begin{nlemma}\label{lemb3}
Let $G$ be a finite group and let $N\unlhd G$. Then
\[k(G)\ \leq\ \frac{k(N)}{|G/N|}\ +\ 2(k(G/N)-1)\max\{|C_{\scl(N)}(g)|\ |\ g\in G-N\}.\]
\end{nlemma}
\begin{beweis}
Let $g_i\in G$ ($i=1,\ldots,k(G/N))$ such that $g_1=1$ and the $\overline{g_i}=g_iN$ are
representatives of the conjugacy classes of $G/N$. Then by \ref{lem1.1}
we have $k(G)=\sum\limits_{i=1}^{k(G/N)}n(C_{G/N}(\overline{g_i}),\Omega_i)$, where $\Omega_i$ is
the set of $N$--orbits on $g_iN$. By the proof of \ref{lem1.1} we have
\[n(C_{G/N}(\overline{g_i}),\Omega_i)\ \leq\ |C_{\scl(N)}(\overline{g_i})|,\]
so that we obtain
\[k(G)\ \leq\ n(G/N,\cl(N))\ +\ \sum_{i=2}^{k(G/N)}|C_{\scl(N)}(\overline{g_i})|.\]
By the Cauchy--Frobenius orbit counting formula we have
\begin{eqnarray*}
n(G/N,\cl(N))&=&\frac{1}{|G/N|}\sum_{gN\in G/N}|C_{\scl(N)}(gN)|\\
              &=&\frac{1}{|G/N|}\sum_{i=1}^{k(G/N)}|G/N:C_{G/N}(\overline{g_i})||C_{\scl(N)}
                 (\overline{g_i})|\\
              &=&\frac{k(N)}{|G/N|}\ +\ \sum_{i=2}^{k(G/N)}\frac{1}{|C_{G/N}(\overline{g_i})|}
                 |C_{\scl(N)}(\overline{g_i})|\\
              &\leq&\frac{k(N)}{|G/N|}\ +\ \sum_{i=2}^{k(G/N)}|C_{\scl(N)}(\overline{g_i})|.
\end{eqnarray*}
Thus altogether
\[k(G)\ \leq\ \frac{k(N)}{|G/N|}\ +\ 2\sum_{i=2}^{k(G/N)}|C_{\scl(N)}(\overline{g_i})|\]
which implies the assertion of the lemma.
\end{beweis}

\section{On the number of fixed conjugacy classes of normal subgroups in certain 
semidirect products}\label{sectionc}\\

The aim of this section is to obtain strong bounds for $|C_{\scl(NV)}(g)|$, where 
$N\unlhd G$, $g\in G$ and $V$ is a faithful $G$--module.\\

We start with an easy lemma.\\

\begin{nlemma}\label{lemc1}
Let $H$ be a finite group, $N\leq G\unlhd H$ and $N\unlhd H$. Let $g\in H$. Then
\[|C_G(g)|\ \leq\ |C_{G/N}(gN)||C_N(g)|.\]
\end{nlemma}
\begin{beweis}
Write $G/N=\{g_iN|i=1,\ldots,|G/N|\}$ for suitable $g_i\in G$. If $h\in C_G(g)$, then
$h=g_ix$ for a unique $i$ and a unique $x\in N$. Now clearly $g_iN=hN\in C_{G/N}(gN)$, so there are
$|C_{G/N}(gN)|$ possibilities for $g_i$. Once $g_i$ is chosen, we see that $g_ix=h=h^g=(g_ix)^g=g_ix_0x^g$ for some
$x_0\in N$ that depends on $g$ and $i$. Hence $x$ is a solution of the equation $[x^{-1},g]=x_0$,
and there are either 0 or exactly $|C_n(g)|$ solutions $x$ for this equation, because if $x_1$ and
$x_2$ are both solutions of the equation, then $x_1^{-1}x_2\in C_N(g^{-1})=C_N(g)$. Hence
the assertion of the lemma follows.
\end{beweis}

\begin{nlemma}\label{lemc2}
Let $L$ be a finite group. Let $H\unlhd L$, and suppose that $|L/H|=p$ is a prime and that 
$H=H_1\times\ldots\times H_p$ for subgroups $H_i$ of $L$ that are permuted by $L/H$, i.e.,
$H_i^g=H_{i+1}$ for $i=1,\ldots,p-1$ and $H_p^g=H_1$, where $L/H=\langle gH\rangle$. Let
$N\leq H$ such that $N^g=N$ and $g^p\in N$, and put $G=\langle N,g\rangle$. Let $N_0= H_1\cap N$
and $N_1=\prod\limits_{i=0}^{p-1}N_0^{g_i}$.
Then obviously $N_0\unlhd N$, $N_1\unlhd G$ and $N_1=\X_{i=0}^{p-1}N_0^{g_i}$. Furthermore, if we put $J=N/N_1$, then
\[|C_{\scl(N)}(g)|\ \leq\ k(J)\cdot k(N_0).\]
\end{nlemma}
\begin{beweis}
By \cite[Lemma 3.1]{guralnick-wan} we have the following elementary formula:
\[|C_{\scl(N)}(g)|\ =\ \frac{1}{|N|}\sum_{n\in N}|C_N(gn)|,\]
and so with \ref{lemc1} we get
\[|C_{\scl(N)}(g)|\ \leq\ \frac{1}{|N|}\sum_{n\in N}|C_J(gnN_1)||C_{N_1}(gn)|.\]
Write $J=\{a_iN_1|i=1,\ldots,|J|\}$ for suitable $a_i\in N$, and also write $M_i=N_0^{g^{i-1}}$ for $i=1,\ldots,p$,
so that $N_1=\X_{i=1}^pM_i$. Then we further have
\begin{eqnarray*}
|C_{\scl(N)}(g)|&\leq&\frac{1}{|N|}\sum\limits_{i=1}^{|J|}\sum\limits_{x_1\in M_1}\ldots\sum\limits_{
                     x_p\in M_p}|C_J(ga_ix_1\ldots x_pN_1)||C_{N_1}(ga_ix_1\ldots x_p)|\\
               &=&\left(\frac{1}{|J|}\sum\limits_{i=1}^{|J|}|C_J(ga_iN_1)|\right)\cdot\left(
                  \frac{1}{|N_1|}\sum\limits_{x_1\in M_1}\ldots\sum\limits_{x_p\in M_p}
                  |C_{N_1}(ga_ix_1\ldots,x_p)|\right).
\end{eqnarray*}
For convenience, call the first factor in the last product $A$ and the second $B$. Then
\[A\ =\ \frac{1}{|J|}\sum_{a\in J}|C_J(ga)|\ =\ \frac{1}{|J|}\sum_{a\in J}|C_{gJ}(a)|,\]
and as it is easy to see that $C_{gJ}(a)$ either is empty or a coset of $C_J(a)$ (see e.g.
\cite[p.176]{gallagher} for the argument), we have
\[A\ \leq\ \frac{1}{|J|}\sum_{a\in J}|C_J(a)|\ =\ k(J).\]
It thus remains to show that $B\leq k(N_0)$.\\

For the moment, fix $i\in\{1,\ldots,|J|\}$ and $x_j\in M_j$ for $j=2,\ldots,p$, and put
$g_0=ga_ix_2\ldots x_p$. Then we clearly have 
\[|C_{N_1}(ga_ix_1x_2\ldots x_p)|=|C_{N_1}(g_0x_1)|,\]
and if we define 
\[U_1=\{z_1\in M_1|z_1z\in C_{N_1}(g_0x_1)\mbox{ for some }z\in M_2\times\ldots\times M_p\}\mbox{, then }
U_1\leq M_1,\] 
and as $g_0x_1$ cyclically
permutes the $M_i$, we see that for each $z_1\in U_1$ there is a unique $z\in M_2\times\ldots\times M_p$
such that $z_1z\in C_{N_1}(g_0x_1)$, so that $|C_{N_1}(g_0x_1)|=|U_1|$. Moreover,
$(g_0x_1)^p\in N$ and clearly $U_1\leq C_{M_1}((g_0x_1)^p)$. Thus altogether
\[|C_{N_1}(ga_ix_1\ldots x_p)|\ \leq\ |C_{M_1}((g_0x_1)^p)|\ =\ |C_{M_1}(g_0^px_1^{g_0^{p-1}}
x_1^{g_0^{p-2}}\ldots x_1^{g_0}x_1)|\ =\ |C_{M_1}(g_0^p x_1)|,\]
where the last equality follows as $x_1^{g_0^j}\in M_{j+1}$ and $M_{j+1}$ centralizes $M_1$
for $j=1,\ldots,p-1$. Moreover we have
\[\sum_{x_1\in M_1}|C_{M_1}(g_0^px_1)|\ =\ \sum_{x_1\in M_1}|C_{g_0^pM_1}(x_1)|\ \leq\ 
\sum_{x_1\in M_1}|C_{M_1}(x_1)|\ =\ |M_1|k(M_1),\]
where the inequality again follows from the fact that $C_{g_0^pM_1}(x_1)$ is either empty or a coset
of $C_{M_1}(x_1)$.\\
With this we finally have
\begin{eqnarray*}
B&\leq&\frac{1}{|N_1|}\sum_{x_2\in M_2}\ldots\sum_{x_p\in M_p}\left(\sum_{x_1\in M_1}
       |C_{M_1}(g_0^px_1)|\right)\\
 &\leq&\frac{1}{|M_1|^p}\sum_{x_2\in M_2}\ldots\sum_{x_p\in M_p}|M_1|k(M_1)\ =\ \frac{1}{|M_1|^p}
       |M_1|^{p-1}|M_1|k(M_1)\\
 &=&k(M_1)\ =\ k(N_0),
\end{eqnarray*}
so that the lemma is proved.
\end{beweis}

We next recall an elementary result essentially obtained in \cite[Lemma 3.3]{kellerkgv2}.
The version that follows, however, has been generalized to include non--coprime actions.\\

\begin{nlemma}\label{leme1}
Suppose that $G$ is a finite group and $V$ is a finite $G$--module. Suppose that
$V=V_1\oplus V_2$ for $G$--modules $V_i$ ($i=1,2$). Let $\lambda_i\in\Irr (V_1)$ ($i=1,\ldots,n(G,V_1)$)
be representatives of the orbits of $G$ on $\Irr (V_1)$. Then
\[k(GV)\ =\ \sum_{i=1}^{n(G,V_1)}k(C_G(\lambda_i)V_2).\]
In particular, $k(GV)\leq n(G,V_1)\cdot\max\{k(C_G(\lambda)V_2\ |\ \lambda\in\Irr (V_1)\}$.\\
Moreover, if $(|G|,|V|)=1$ and $v_i\in V_1$ ($i=1,\ldots,n(G,V_1)$) are representatives of the
orbits of $G$ on $V_1$, then
\[k(GV)=\sum_{i=1}^{n(G,V_1)}k(C_G(v_i)V_2)\leq n(G,V_1)\max\{k(C_G(v)V_2)\ |\ v\in V_1\}.\]
\end{nlemma}
\begin{beweis}
First observe that $n(G,V_1)=n(G,\Irr (V_1))$ and $n(G,V)=n(G,\Irr (V))$ by Brauer's permutation
lemma. Note that any $\lambda\in\Irr(V)$ can be extended to its inertia group in $GV$,
and therefore if $\mu_i$ ($i=1,\ldots,n(G,V)$) are representatives of the orbits of $G$ on
$\Irr (V)$, then with Gallagher's result \cite[Corollary (6.12)]{isaacs} we conclude that
\begin{eqnarray*}
k(GV)&=&|\Irr (GV)|\ =\ \sum_{i=1}^{n(G,V)}k(C_G(\mu_i))\\
     &=&\frac{1}{|G|}\sum_{\mu\in\Irr (V)}|C_G(\mu)|k(C_G(\mu))\\
     &=&\frac{1}{|G|}\sum_{\mu_1\in\Irr (V_1),\mu_2\in\Irr (V_2)}|C_G(\mu_1\mu_2)|k(C_G(\mu_1\mu_2))\\
     &=&\sum_{\mu_1\in\Irr (V_1)}\frac{|C_G(\mu_1)|}{|G|}\left(\frac{1}{|C_G(\mu_1)|}\sum_{
        \mu_2\in\Irr (V_2)}|C_{C_G(\mu_1)}(\mu_2)|k(C_{C_G(\mu_1)}(\mu_2))\right)\\
     &=&\sum_{\mu_1\in\Irr (V_1)}\frac{|C_G(\mu_1)|}{|G|} k(C_G(\mu_1)V_2)\\
     &=&\sum_{i=1}^{n(G,V_1)}k(C_G(\lambda_i)V_2),
\end{eqnarray*}
as wanted. The remaining statements are immediate and well--known consequences of 
the first one.
\end{beweis}

We also need a result on the number of orbits.\\

\begin{nlemma}\label{lema3}
Let $G$ be a finite group and let $V$ be a finite $G$--module. Let $N\unlhd G$. Then
\[n(G,V)\ \leq\ \left(\frac{k(GV)|V|}{k(G)}\right)^\frac{1}{2}.\]
\end{nlemma}
\begin{beweis}
First note that $n(G,\Irr (V))=n(G,V)$ by Brauer's permutation lemma.
As every $\lambda\in\Irr (V)$ can be extended to its inertia group in $GV$, we have, if
the $\lambda_i$ are representatives of the orbits of $G$ on $\Irr (V)$, that
\begin{eqnarray*}
k(GV)&=&\sum_{i=1}^{n(G,V)}k(C_G(\lambda_i))\\
     &\geq&\sum_{i=1}^{n(G,V)}\frac{k(G)}{|G:C_G(\lambda_i)|}\ =\ k(G)n(G,V)
           \frac{1}{n(G,V)}\sum_{i=1}^{n(G,V)}\frac{1}{|G:C_G(\lambda_i)|}\\
     &\geq&k(G)n(G,V)^2\left(\sum_{i=1}^{n(G,V)}|G:C_G(\lambda_i)|\right)^{-1}
           =k(G)n(G,V)^2|V|^{-1}, 
\end{eqnarray*}
where the first inequality follows from Ernest's result (see \cite[Problem E17.2]{huppertcharacters})
and the second inequality follows from the arithmetic--harmonic--means inequality. Therefore
$n(G,V)\leq\left(\frac{k(GV)|V|}{k(G)}\right)^{1/2}$, and so we are done.
\end{beweis}

We now can prove an important auxiliary result.\\

\begin{nlemma}\label{lemc4}
Let $G$ be a finite group and $V$ be a finite faithful $G$--module. Suppose that $p$ is a prime
and $V=V_1\oplus\ldots\oplus V_p$ for subspaces $V_i$ which are permuted nontrivially by $G$,
and put $N=\bigcap\limits_{i=1}^pN_G(V_i)\unlhd G$. Moreover, assume that $G/N=\langle gN\rangle$ 
is cyclic of order $p$. Put 
\[K_i=C_N\left(\bigoplus_{j=1;j\not=i}^pV_j\right)\unlhd N\] 
for $i=1,\ldots,p$, so then $N_1=K_1\ldots K_p=K_1\times\ldots\times K_p\unlhd G$.\\
Put  $N_0=K_1$, and $J=N/N_1$. Then 
\[|C_{\scl(NV)}(g)|\leq k(J)\cdot k(N_0V_1).\] 
Put $W_2=V_2\oplus\ldots\oplus V_p$, $U_1=N_G(V_1)/C_G(V_1)$ and
\[M=\max\{k((C_N(\lambda_1)N_0/N_0)W_2\ |\ \lambda_1\in\Irr(V_1)\}\] 
and $m=\max\{k(T)|T\leq N_0\}$.
Then for any $S\leq U_1$ with $k(S)=\max\{k(U)|U\leq U_1\}$ we have 
\[k(NV)\ \leq\ \left(\frac{k(SV_1)|V_1|}{k(N/N_0)^\frac{1}{p-1}}\right)^{1/2}\cdot M\cdot m\ \leq\ 
\left(\frac{k(SV_1)|V_1|}{k(J)^\frac{1}{p-1}}\right)^{1/2}\cdot M\cdot m\]
\end{nlemma}
\begin{beweis}
Note that if we put $H_i=(N_G(V_i)/C_G(V_i))V_i$ for $i=1,\ldots,p$ and $H=H_1\times\ldots\times H_p$
(so that $NV\ug H$)
and $L=\langle H,g\rangle$, then (after possibly relabeling the $H_i$) the hypotheses of
\ref{lemc2} are fulfilled, so it follows easily that  
$|C_{\scl(NV)}(g)|\leq k(J)\cdot k(N_0V_1)$ which proves the first inequality 
that we have to establish.\\
To prove the second one, put $W_1=V_1$ and
observe that $W_1$ and $W_2$ are $N$--modules. By \ref{leme1},
we have 
\[k(NV)\ \leq\ n(N,W_1)\cdot M_1,\mbox{ where }M_1=\max\{k(C_N(\lambda_1)W_2)\ |\ \lambda_1\in\Irr(V_1)\}.\]
Now let $S\leq U_1$ such that $k(S)=\max\{k(U)|U\leq U_1\}$. Then by \ref{lema3} we have
\[n(N,W_1)\ =\ n(N/C_N(W_1),V_1)\ \leq\ n(U_1,V_1)\ \leq\ n(S,V_1)\ \leq\ \left(\frac{k(SV_1)|V_1|}{
k(S)}\right)^\frac{1}{2}.\]
Recall that $C_N(W_2)=N_0$.
Furthermore it is easy to see that with $X_i=\bigoplus\limits_{j=2}^iV_j$ ($i=2,\ldots,p$) we have
\[k(N/N_0)\ =\ k(N/C_N(W_2))\ \leq\ \prod\limits_{i=2}^pk(C_{N/N_0}(X_{i-1})/C_{N/N_0}(X_i))\ 
\leq\ k(S)^{p-1},\]
where the second inequality follows by the choice of $S$. Thus $k(S)\geq (k(N/N_0))^\frac{1}{p-1}$
and hence
\[n(N,W_1)\ \leq\ \left(\frac{k(SV_1)|V_1|}{k(N/N_0)^{1/p}}\right)^\frac{1}{2}\ \leq\ 
\left(\frac{k(SV_1)|V_1|}{k(J)^\frac{1}{p-1}}\right)^\frac{1}{2}.\]
To complete the proof of the lemma, it remains to show that $M_1\leq M\cdot m$.
For any subgroup $T\leq N$ we have
\[k(TW_2)\ \leq\ k((T/C_T(W_2))W_2)\cdot k(C_T(W_2))\ \leq\ k((TN_0/N_0)W_2)\cdot m,\]
and so the assertion follows and we are done.
\end{beweis}

\section{The coprime case}\label{sectiond}

In this section we study what our results yield in the situation of the classical $k(GV)$--problem. For this,
first recall the following result by Gambini and Gambini--Weigel, as stated in 
\cite[Theorem 2.1]{gluck-magaard}.\\

\begin{nsatz}\label{theoremd1}
Let $G$ be a finite group and $W$a faithful primitive finite $G$--module with $(|G|,|W|)=1$.
Then 
\[|G|\leq |W|\log_2|W|,\]
except when $|W|=7^4$ and $G$ is ${\rm Sp}(4,3)$ or $Z_3\times {\rm Sp}(4,3)$.
\end{nsatz}

\begin{nlemma}\label{lemd2}
Let $G$ be a finite group and let $V$ be a finite faithful $G$--module with $(|G|,|V|)=1$,
and assume that $k(UX)\leq |X|$ whenenver $U$ is a finite group, $X$ is a faithful $U$--module
with $(|U|,|X|)=1$ and $|UX|<|GV|$. Suppose that $W<V$ and $H<G$ are such that $H=N_G(W)$, $W$
is primitive as $H$--module, and $V=W^G$ is induced from $W$. So we can write $V=V_1\oplus\ldots\oplus V_n$
for some $n>1$ and subspaces $V_i$ that are permuted faithfully by $G/N$, where $N=\bigcap\limits_{
g\in G}H^g$ and $V_1=W$. 
Let $p$ be a prime and
$g\in G-N$ such that $g^p\in N$. Let $f$ be the number of $p$--cycles in the permutation action
of $g$ on $\{V_1,\ldots,V_n\}$, so that $g$ normalizes $n-pf$ of the $V_i$. Clearly $f\geq 1$. Put
\[B\ =\ \left\{\begin{array}{ll}
6&\mbox{if }|V_1|=7^4\mbox{ and }N/C_N(V_1)\mbox{ is isomorphic to Sp}(4,3)\mbox{ or }Z_3\times {\rm Sp}(4,3)\\
1&\mbox{otherwise}
\end{array}\right.\]
Then
\[|C_{\scl(NV)}(g)|\ \leq\ A_i^f|V_1|^{n-pf},\ i=1,2,\]
for each of the following $A_i$:
\[\begin{array}{l}
A_1=B|V_1|^2\log_2|V_1|\mbox{ and}\\
A_2=|V_1|^\frac{2p^2+1}{2p+1}(B\log_2|V_1|)^\frac{2p-2}{2p+1}
\end{array}\]
Note that for $p=2$, $A_2=|V_1|^\frac{9}{5}(B\log_2|V_1|)^\frac{2}{5}$.
\end{nlemma}
\begin{beweis}
Note that if we put $D_i=(N/C_N(V_i))V_i$, then $NV\ug M:=\X_{i=1}^nD_i$ with $G/N$ permuting the factors
of this direct product. 
Now relabel the $D_i$ such that ${\cal O}_i=\{D_{(i-1)p+1},\ldots,D_{ip}\}$ ($i=1,\ldots,f$)
are the orbits of $\langle g\rangle$ on $\{D_1,\ldots,D_n\}$ of size $p$, and
${\cal O}_i=\{D_{(p-1)f+i}\}$ ($i=f+1,\ldots,n-(p-1)f$) are the remaining orbits.
Put $M_i=\X_{j\mbox{ \scriptsize with }D_j\in{\cal O}_j}D_j$ for $i=1,\ldots,n-(p-1)f$ and consider
the group $G_0=\langle g,M\rangle$.\\
This clearly satisfies the hypothesis of \ref{lem1.2}, and if we define the $k_i$ as in
\ref{lem1.2}(b), then by our hypothesis $k_i\leq |V_i|=|V_1|$ for $i=f+1,\ldots,n-(p-1)f$,
and so by \ref{lem1.2}(b) we have
\[|C_{\scl(NV)}(g)|\leq |V_1|^{n-pf}\prod_{i=1}^fk_i.\] 
Thus it remains to show that $k_i\leq A_1$
and $k_i\leq A_2$ for $i=1,\ldots,f$. For this we clearly may assume that $f=1$, $n=p$ and $G$
is embedded in 
$H=\langle g,H_1\times\ldots\times H_p\rangle$, where $H_i=N_G(V_i)/C_G(V_i)$ for
$i=1,\ldots,p$ and $H_i^g=H_{i+1}$ for $i=1,\ldots,p-1$ and $H_p^g=H_1$, and $N\leq H$, and
we have to show that for $C:=C_{\scl(NV)}(g)$ we have
\[|C|\ \leq \ A_i\mbox{ for }i=1,2.\]
Now by \ref{theoremd1} we have $|D_1|\leq B|V_1|^2\log_2|V_1|$ and so by \ref{lem2}
we have
\[|C|\ \leq\ B|V_1|^2\log_2|V_1|.\]
This gives the first part of the lemma.\\
Next let $N_0$, $N_1$ and $J$ be as in \ref{lemc4}. With \ref{lemc4} and our hypothesis we obtain
that
\[(1)\quad |C|\leq k(J)|V_1|\]
and 
\begin{eqnarray*}
|C|\leq k(NV)&\leq&\frac{|V_1|}{k(J)^\frac{1}{2(p-1)}}\cdot |V_1|^{p-1}\cdot\max\{k(T)|T\leq N_0\}\\
             &\leq&\frac{|V_1|^p}{k(J)^\frac{1}{2(p-1)}}\cdot |N_0|.
\end{eqnarray*}
Observe that as $C_N(V_2\oplus\ldots\oplus V_p)=N_0$, by \ref{theoremd1} we have
\[|J||N_0|^{p-1}=|N/N_0|=|N/C_N(V_2\oplus\ldots\oplus V_p)|\leq |H_1|^{p-1}\leq (B|V_1|\log_2|V_1|)^{p-1}\]
and thus
\[|N_0|\ \leq\ \frac{B|V_1|\log_2|V_1|}{|J|^\frac{1}{p-1}},\]
so that we further get
\[(2)\quad |C|\ \leq\ \frac{B|V_1|^{p+1}\log_2(|V_1|)}{k(J)^\frac{1}{2(p-1)}|J|^\frac{1}{p-1}}\
\leq\ \frac{B|V_1|^{p+1}\log_2(|V_1|)}{k(J)^\frac{3}{2(p-1)}}.\]
Now the upper bounds in (1) and (2) are equal if and only if 
\[(3)\quad k(J)\ =\ \left(B|V_1|^p\log_2|V_1|\right)^\frac{2p-2}{2p+1}.\]
Therefore either (1) or (2) will always yield a bound less than or equal to the one we obtain
in case that $k(J)$ has the critical value in (3), therefore we always have
\begin{eqnarray*}
(4)\quad |C|&\leq&(B|V_1|^p\log_2|V_1|)^{\frac{2p-2}{2p+1}}\cdot |V_1|\\
            &=&|V_1|^{\frac{2p^2+1}{2p+1}}(B\log_2|V_1|)^{\frac{2p-2}{2p+1}}
\end{eqnarray*}
So the lemma is proved.
\end{beweis}

Note that the $A_1$--bound in the previous lemma, which was relatively easy to establish, is
always much stronger than the $A_2$--bound, with the only exception of $p=2$, where the
$A_1$--bound is trivial and useless, while the $A_2$--bound is nontrivial, albeit quite weak.\\
In view of the applications of \ref{lemd2}, it would be highly desirable to improve the bound
for $p=2$; the current $A_2$--bound seems to be much too large, and in fact something like
\[|C|\ \leq\ |V_1|^\frac{2}{3}\]
instead of (4) should be possible.
While our general bounds are larger than necessary, in specific situations, when more
detailed information on the groups is available, such as good bounds for $k(J)$, then
the formulas (1), (2) in the proof of \ref{lemd2} will yield much better results, as we shall see in Section \ref{sectione}. 
This is already
so in case that $N/C_N(W_1)$ is isomorphic to ${\rm Sp}(4,3)$ or $Z_3\times {\rm Sp}(4,3)$, so
that better bounds than the ones in the previous lemma can be obtained in that case, although
we will not pursue this further here.\\

Next we look at the most general reduction of the imprimitive case of the $k(GV)$--problem
that we can get here.\\

\begin{nsatz}\label{lemd4}
Let $G$ be a finite group and $V$ be a finite faithful $G$--module with $(|G|,|V|)=1$. Assume
that $k(UX)\leq |X|$ whenever $U$ is a finite group, $X$ is a faithful $U$--module with
$(|U|,|X|)=1$ and $|UX|<|GV|$. Suppose further that $W<V$ and $H<G$ are such that $H=N_G(W)$,
$W$ is primitive as $H$--module, and $V=W^G$ is induced from $W$. Put $\overline{H}=H/C_H(W)$. Then the
following hold:\\
(a) If $|W|\geq 2^{47}$ then $k(GV)\leq |V|$.\\
(b) If $k(HV)\leq |V|-(3^{(n-1)/2}+1)|V|^\frac{9}{10}(6\log_2|W|)^\frac{n}{5}$, then $k(GV)\leq |V|$.
Moreover, if $k(\overline{H}W)\leq\frac{|W|}{2}$ and $|W|\geq 2^{19}$, then $k(GV)\leq |V|$.
\end{nsatz}
\begin{beweis}
Let $N=\bigcap\limits_{g\in G}H^g\unlhd G$. Then we can write $V=V_1\oplus\ldots\oplus V_n$
for $n=|G:H|$ and submodules $V_i\leq V$ such that $V_1=W$, and $G/N$ permutes the $V_i$
transitively and faithfully.\\

(a) If $g\in G$ has at most $\frac{n}{2}$ fixed points
in its permutation action on $\{V_1,\ldots,V_n\}$, then by \ref{lemd2} we know that
with $B$ as in \ref{lemd2} we have
\[|C_{\scl(NV)}(g)|\ \leq\ \left(|W|^\frac{9}{5}(B\log_2|W|)^\frac{2}{5}\right)^\frac{n}{4}|W|^\frac{n}{2}\ 
=\ |W|^{\frac{19}{20}n}(B\log_2|W|)^\frac{n}{10}\ =:\ C\]
(as we clearly may assume that $g$ is of prime order when checking this).\\
Now by \ref{lemb2} (with $\delta=1$ and $f(x)=x$) we are done if
\[k(UN/N)\ \leq\ \frac{1}{\sqrt{n+1}}\left(\frac{|V|}{C}\right)^\frac{1}{2}\mbox{ for all }
U\leq G,\]
and as for $n\ne 2$ we have $k(UN/N)\leq 3^{(n-1)/2}$ for all $U\leq G$, as $UN/N$ is isomorphic to a subgroup of $S_n$,
it suffices to have
\[3^{\frac{n-1}{2}}\ \leq\ \frac{1}{\sqrt{n+1}}\left(\frac{|V|}{C}\right)^\frac{1}{2}\ =\ \frac{1}{\sqrt{n+1}}
\frac{|W|^\frac{n}{40}}{(B\log_2|W|)^\frac{n}{20}}\]
and also (in case that $n=2$) that
\[2\ \leq\ \frac{1}{\sqrt{3}}\ \frac{|W|^\frac{1}{20}}{(B\log_2|W|)^\frac{1}{10}}.\]
This is the case for $|W|\geq 2^{47}$ (as $B=1$ in this case), as can easily be verified, so (a) is proved.\\

(b) If $g\in G$ permutes the $V_i$ ($i=1,\ldots,n$) fixed point freely, then by \ref{lemd2} we know that
\[|C_{\scl(NV)}(g)|\leq\left(|W|^\frac{9}{5}(B\log_2|W|)^\frac{2}{5}\right)^\frac{n}{2}=:D,\]
where $B$ is as in \ref{lemd2}.\\
By \ref{lemb1} we have
\begin{eqnarray*}
k(GV)&\leq&k(HV)+k(G/N)\max\{|C_{\scl(NV)}(g)|\ |\ g\in G-\bigcup_{x\in G}H^x\}\\
     &\leq&k(HV)+k(G/N)\max\{|C_{\scl(NV)}(g)|\ |\ g\in G\mbox{ permutes the }V_i\mbox{ fixed point freely.}\}
\end{eqnarray*}
As $G/N\ug S_n$, again by \cite{maroti} we have $k(G/N)\leq \left\lceil3^{(n-1)/2}\right\rceil\leq 2^{n-1}$ (where
$\lceil x\rceil$ denotes the upper integer part of $x$), and so
we conclude that
\begin{eqnarray*}
k(GV)&\leq&k(HV)+\left\lceil 3^{(n-1)/2}\right\rceil\  D\\
     &=&k(HV)+\left\lceil 3^{(n-1)/2}\right\rceil\ |W|^{\frac{9}{10}n}(B\log_2|W|)^\frac{n}{5},
\end{eqnarray*}
and so by our hypothesis the first assertion of (b) follows.\\
To prove the second one, first note that if $k(\overline{H}W)\leq\frac{|W|}{2}$, then
\begin{eqnarray*}
k(HV)&\leq&k(\overline{H}V_1)\cdot k(C_H(V_1)(V_2\oplus\ldots\oplus V_n))\leq\frac{|V_1|}{2}
           \cdot |V_2\oplus\ldots\oplus V_n|\\
     &=&\frac{|V|}{2},
\end{eqnarray*}
and as $|W|>7^4$, clearly $B=1$ here, so we obtain
\[k(GV)\leq\frac{|V|}{2}+2^{n-1}|W|^{\frac{9}{10}n}(\log_2|W|)^\frac{n}{5}.\]
Thus $k(GV)\leq |V|$ if 
\[2^{n-1}|W|^{\frac{9}{10}n}(\log_2|W|)^\frac{n}{5}\ \leq\ \frac{|W|^n}{2}\]
which is equivalent to 
\[2^{10}(\log_2|W|)^2\ \leq\ |W|,\] 
and this holds for $|W|\geq 2^{19}$. So the
theorem is proved.
\end{beweis}

So this is a general reduction of the imprimitive case of the 
$k(GV)$--problem to ``small" cases.
For large primes $p=\char (V)$, this result even provides a complete reduction of the
imprimitive case to the primitive case, saying that a minimal counterexample to the $k(GV)$--problem
must be primitive. This is an improvement of the corresponding part in the proof of 
\cite[Theorem 4.1]{kellerkgv2}.
(For ways to treat the primitive case for large primes, see \cite{kellerkgv2}.)\\
It would be nice if one could refine the methods here, in particular improve the bounds in
\ref{lemd2}, so as to further reduce the $2^{47}$ in \ref{lemd4} and reach a general reduction
of the $k(GV)$--problem to primitive actions.\\
It would also be interesting to know whether with methods as the ones employed here it is possible
(at least for large $p$) to reduce the problem further to tensorprimitive modules $V$.\\

\section{The last case of the $k(GV)$--problem}\label{sectione}\\

While \ref{lemd4} seems to imply that our techniques only work for large primes $p=\char (V)$,
we will now see that they are also quite powerful in ``small'' situations. We demonstrate
this by providing a new, short proof of the $k(GV)$--problem in the situation that turned out to
be the most difficult in the original proof of the $k(GV)$--problem and that occupied all of 
\cite{gmrs}. Here $p=5$ and $V$ is induced from the irreducible module of order $5^2$ of a
5--complement of $\GL (2,5)$.\\

\begin{nlemma}\label{leme2}
Let $G$ be a finite group and $V$ be a finite faithful $G$--module. Let $p$ be a prime, and
suppose that $V=V_1\oplus\ldots\oplus V_p$ for subspaces $V_i$ which are permuted nontrivially by
$G$. Assume that $G/N=\langle gN\rangle$ is cyclic of order $p$, where $N=\bigcap\limits_{i=1}^pN_G(V_i)\unlhd G$.
Suppose further that $|V_1|=5^2$ and let $L$ be a 5--complement of $\GL (2,5)$. Assume that
$U_1:=N_G(V_1)/C_G(V_1)$ is isomorphic to a subgroup of $L$ in its natural action on $V_1$. Then
\[|C_{\scl(NV)}(g)|\ \leq\ |V|^{0.74}.\]
\end{nlemma}
\begin{beweis}
Put $C=C_{\scl(NV)}(g)$. By \ref{lem2} we have 
\[|C|\ \leq\ |LV_1|\ =\ 96\cdot 25\ =\ 2400.\]
Hence $|C|\leq |V|^{0.74}$ if $p\geq 5$, as can easily be checked.\\

Let $p\leq 3$. Let $N_0\leq N$ and $N_1\unlhd G$ be as in \ref{lemc4}, and write
$\overline{N}=N/C_N(V_1)$ and observe that we may consider $N_0$ to be a subgroup of
$\overline{N_1}$ as $N_0$ acts faithfully on $V_1$.
Put $J=N/N_1$.\\

Let $p=3$. If $N/C_N(V_1)$ is not isomorphic to $L$, then 
\[|N/C_N(V_1)|\leq\frac{96}{2}=48,\]
and then as for the primes $\geq 5$ by \ref{lem2} we conclude that 
\[|C|\leq48\cdot 25\leq|V|^{0.74}=25^{2.22}.\]
So we may assume that $\overline{N}\cong L$.\\
Now if $|N_0|\geq 8$,
then $N_0$ contains $\overline{G_1}''\cong Q_8$ (the quaternion group of order 8), 
and then it is easy to see that $k(J)\leq 50$,
so by \ref{lemc4} we have 
\[|C|\leq k(J)k(N_0V_1)\leq 50\cdot 25\leq |V|^{0.74},\]
so that we are
done in this case. Hence $|N_0|\in\{1,2,4\}$.\\
If $|N_0|=4$, then $\overline{N}/N_0\cong S_4$ and
thus again $k(J)\leq k(S_4)^2=5\cdot 5=25$, so again by \ref{lemc4} we have 
\[|C|\leq 25\cdot 25\leq |V|^{0.74}.\]
In the remaining cases we use \ref{leme1}. Note that as $\overline{N}\cong L$, we have
$n(N,V_1)=2$, so if $|N_0|=1$, then $N$ acts faithfully on $V_2+V_3$, and by \ref{leme1} for any
$0\not=v_1\in V_1$ we have 
\begin{eqnarray*}
|C|&\leq&k(NV)\ =\ k(N(V_2\oplus V_3))\ +\ k(C_N(v_1)(V_2\oplus V_3))\\
   &\leq&k(N/C_N(V_2)V_2)k(C_N(V_2)V_3)\ +\ 25^2\ \leq\ 20\cdot 25\ +\ 25^2\ =\ 1125\ \leq\ |V|^{0.74},
\end{eqnarray*}
as wanted. Thus let $|N_0|=2$. Let $0\not=v_1\in V_1$.
Then $C_N(v_1)$ acts faithfully on $V_2\oplus V_3$, and    
$|C_N(V_2\oplus V_3)|=2$, and thus
\[J=N/N_1=N/(N_1N_0)\cong (N/N_0)/(N_1/N_0),\]
and as
\[N/N_0\ug N/C_N(V_2)\times N/C_N(V_3)\ug L\times L\]
and $|N_1/N_0|=|N_0|^2=4$, we see that
\[J\ug L/Z(L)\times L/Z(L)\] and
so $|J|\ |\ 48^2$. If $|J|\ |\ 
\frac{48^2}{2}$, then it is clear from the structure of $L/Z(L)\cong S_4\times C_2$ that
$k(J)\leq 50$, and then as $k(N_0V_1)=14$, by \ref{lemc4} we have $|C|\leq 14\cdot 50=700\leq |V|^{0.74}$. 
Hence we may assume that $J\cong L/Z(L)\times L/Z(L)$. Then
\[k(N(V_2\oplus V_3))\ \leq\ k((N/C_N(V_2\oplus V_3))(V_2\oplus V_3))\cdot k(C_N(V_2\oplus V_3))\
\leq\ 20\cdot 20\cdot 2\ =\ 800.\]
Moreover, $C_N(v_1)/C_N(V_1)\cong C_4$, and $C_N(V_1)/C_{N_1}(V_1)\cong S_4\times C_2$, and
$|C_N(V_1)\cap C_N(V_i)|=2$ for $i=2,3$, and so $C_N(V_1)/C_{C_N(V_1)}(V_2)\cong L$. Hence
$n(C_N(V_1),V_2)=2$, and by \ref{leme1} we conclude that 
\begin{eqnarray*}
k(C_N(v_1)(V_2\oplus V_3))&\leq&4\cdot k(C_N(V_1)(V_2\oplus V_3))\\
                          &\leq&4\cdot 2\cdot\max\{k(C_{C_N(V_1)}(v_2)V_3)\ |\ v_2\in V_2\}\\
                          &\leq&4\cdot 2\cdot 25\cdot 2\ =\ 400.
\end{eqnarray*}
Thus altogether by \ref{leme1} we have
\begin{eqnarray*}
|C|\ \leq\ k(NV)&=&k(N(V_2\oplus V_3))\ +\ k(C_N(v_1)(V_2\oplus V_3))\\
                &\leq&800\ +\ 400\ =\ 1200\ \leq\ |V|^{0.74}
\end{eqnarray*}
which concludes the case $p=3$.\\

It remains to consider the case $p=2$. Here we have to show that $|C|\leq 117$.\\
Now if $k(J)\leq 4$, then again by \ref{lemc4} we have $|C|\leq 4\cdot 25=100$ and we are done.
Thus from now on let $k(J)\geq 5$.\\
If $3|\ |N_0|$, then $J$ is a 2--group and thus $|J|\geq 8$, so $L$ has a section of order
24 with a normal Sylow 3--subgroup, which contradicts the structure of $L$. Thus $3\not| |N_0|$.\\
Next suppose that $3\not| |J|$, so $3\not| |U_1|$ and $U_1$ is a 2--group, more precisely a subgroup
of $S:=C_4\wr C_2$ (which is a 2--Sylow subgroup of $L$). Observe that $S'$ is of order 4 and acts
fixed point freely on $V_1$.\\
Assume that $Z:=Z(F(L))\leq N_0$ (up to isomorphism). Then $|N_0|\in\{2,4\}$ (as $k(J)\geq 5$).
If $|N_0|=4$, it is easy to check that then $k(N_0V_1)\leq 16$, so by \ref{lemc4} we have 
$|C|\leq 16k(J)$ which forces $k(J)=8$ (otherwise we are done), and so $J$ is abelian of order 8. But then $N_0=S'$ and thus
$k(N_0V_1)=10$, so that \ref{lemc4} yields the wanted conclusion here. If $|N_0|=2$, then
$k(N_0V_1)=14$, so by \ref{lemc4} we may assume that $k(J)>8$ forcing $|J|=16$. Now we use
\ref{leme1} which in our coprime situation means that if $v_1=0,v_2,v_3$ are representatives of the
three orbits of $U_1$ on $V_1$, then
\[k(NV)\ =\ k(NV_2)+k(C_N(v_2)V_2)+k(C_N(v_3)V_2),\]
and we may assume that $C_{U_1}(v_2)$ is cyclic of order 4 and $C_{U_1}(v_3)$ is of order 2.
Hence $C_N(v_2)\cong C_4\times C_2$ and $C_N(v_3)$ is elementary abelian of order 4, and both $C_N(v_2)$
and $C_N(v_3)$ act faithfully on $V_2$ (because if, say, $C_N(v_2)$ contained an element $x$
acting trivially on $V_2$, then $x\in N_0$, contradicting the fact that $N_0$ acts fixed point
freely on $V_1$),
and both $C_N(v_2)$ and $C_N(v_3)$
contain an involution acting fixed point freely on
$V_2$. Thus it is easy to check that $k(C_N(v_2)V_2)\leq 25$ and $k(C_N(v_3)V_2)\leq 25$, and as
\[k(NV_2)\leq 2\cdot k((N/C_N(V_2))V_2)\leq 2\cdot 25=50,\]
altogether we have $|C|\leq k(NV)\leq 25+25+50=100$,
as wanted.\\ 
Therefore to complete the case $3\not| |J|$ we may assume that $Z\not\leq N_0$.
But from the structure of $L$ it is then clear that then $|J|\leq 4$, a contradiction.\\

Hence for
the rest of the proof we may assume that $3|\ |J|$. From the structure of $L$ and since $k(J)\geq 5$
it then follows that $|N_0|\leq 8$.\\
If $|N_0|=8$, then necessarily $N_0\cong Q_8$ acts fixed point
freely on $V_1$ and thus $k(N_0V_1)=8$, so by \ref{lemc4} we have $|C|\leq k(J)k(N_0V_1)\leq
|J|\cdot 8=12\cdot 8=96$, so we are done here.\\
If $|N_0|=4$, then $N_0$ is cyclic of order 4 and
acts fixed point freely on $V_1$, so that $k(N_0V_1)=10$. Moreover $k(J)\leq 10$ and so again
by \ref{lemc4} we are done.\\
Next let $|N_0|=2$. Then $N_0$ acts fixed point freely on $V_1$, so that
$k(N_0V_1)=14$. Thus by \ref{lemc4} we may assume that $k(J)\geq 9$, which implies that $U_1\cong L$
is as large as possible. Then by \ref{leme1} we have 
\[k(NV)=k(NV_2)+k(C_N(v)V_2)\] 
for any
$0\not=v\in V_1$, as $U_1$ has only one nontrivial orbit on $V_1$. As clearly $k(NV_2)$ and
$k(C_N(v)V_2)$ are bounded above by $2\cdot |V_2|=50$, we are done here as well.\\
So finally
let $|N_0|=1$. Hence $N$ acts faithfully on $V_1$ and on $V_2$, and as $k(UV_1)\leq |V_1|=25$ for any
$U\leq N$, by \ref{leme1} we have $k(NV)\leq 25n(N,V_1)$, so that we are done whenever $n(N,V_1)\leq 4$.
Thus let $n(N,V_1)\geq 5$. Then from the structure of $H$ and its action on $V_1$ it is clear that
$N$ must be cyclic of order 3 or 6 and act fixed point freely on $V_1$. As $J\cong N$ here
and $k(J)\geq 5$, only the case $|N|=6$ remains, and then $|C|\leq k(NV)=110$. This completes the proof
of the lemma.
\end{beweis}

Now we can prove the main result of this section. Observe that this includes the main result of
\cite{gmrs}, which constituted the last and in some sense most difficult case of the $k(GV)$--problem.\\

\begin{nsatz}\label{theoreme3}
Let $G$ be a finite $5'$--group and $V$ be a faithful GF$(5)$--module such that $V$ is induced
from a $G_1$--module $W$, where $G_1$ is a suitable subgroup of $G$, $|W|=25$ and
$G_1/C_{G_1}(W)\not=1$ is isomorphic to a subnormal subgroup of $L$, where $L$ is a 5--complement
in $\GL (2,5)$. Suppose that whenever $U\leq G$ and $X\leq V$ is a $U$--module with
$|UX|<|GV|$, then $k(UX)\leq|X|$.
Then
\[k(GV)\leq|V|.\]
\end{nsatz}
\begin{beweis}
Put $n=|G:G_1|$. Clearly we may assume that $n>1$. Then $V=V_1\oplus\ldots\oplus V_n$ for subspaces
$V_i\cong W$ that are permuted transitively by $G$. Write $H=N_G(V_1)$ and 
$N=\bigcap\limits_{x\in G}H^x\unlhd G$. Now let $g\in G-\bigcup\limits_{x\in G}H^x$ be of prime 
order, so that in particular $g$ has no fixed point in its permutation action on $\{V_1,\ldots,V_n\}$.
Then applying \ref{lem1.2}(b) and \ref{leme2} to the group $\langle g,NV\rangle$ yields
\[|C_{\scl(NV)}(g)|\ \leq\ |V|^{0.74}\quad (1)\]
and hence (1) holds for all $g\in G-\bigcup\limits_{x\in G}H^x$. \\

First suppose that $n=2$. Then $|G/N|=2$, and \ref{lem1.1}, together with (1), yields
\[k(GV)\ \leq\ k(NV)+|V|^{0.74}\quad (2).\]
Put $N_1=N/C_N(V_1)$, then 
\[k(NV)\leq k(N_1V_1)\cdot k(C_N(V_1)V_2)\leq k(N_1V_1)\cdot 25,\] 
the second inequality following by our hypothesis. So (2) yields $k(GV)\leq|V|=625$ unless
$k(N_1V_1)\geq 21$, but it can easily be checked that this happens only when $V_1$ is
reducible as $N_1$--module and either $N_1=1$ or $N_1$ is cyclic of order 4. But as $N=G_1$ here,
we have $N_1\cong G_1/C_{G_1}(W)\not=1$ by hypothesis, and if $|N_1|=4$ and $V_1$
is reducible as $N_1$--module, then $N_1$ is not subnormal in (its copy in) $L$, again contradicting
our hypothesis.\\

Thus for the rest of the proof we may assume that $n\geq 3$.
Now by \ref{lemb1} and (1) we have 
\begin{eqnarray*}
k(GV)&\leq&k(HV)+k(G/N)\max\{|C_{\scl(NV)}(g)|\ \big|\ g\in G-\bigcup\limits_{x\in G}H^x\}\\
     &\leq&k(HV)+k(G/N)|V|^{0.74}\quad (3).
\end{eqnarray*}
Recall that if $S\leq S_n$, then $k(S)\leq (\sqrt{3})^{n-1}$
for $n\geq 3$ (see \cite{maroti}). With this, (3) becomes
\[k(GV)\leq k(HV)+3^{\frac{1}{2}(n-1)}|V|^{0.74}\quad (4)\]
Next we have to bound $k(HV)$. Note that $V$ is reducible as $H$--module, and so with
$H_1:=H/C_H(V_1)$ we have 
\[k(HV)\leq k(H_1V_1)\cdot k(C_H(V_1)(V_2\oplus\ldots\oplus V_n))\leq k(H_1V_1)\cdot 25^{n-1}\]
(again by our hypothesis).
Now in the case that $k(H_1V_1)\geq 21$ as for $n=2$ we run into a contradiction, so we
may assume that $k(H_1V_1)\leq 20$. Thus from (4) we get
\[k(GV)\leq\frac{4}{5}|V|+3^{\frac{1}{2}(n-1)}|V|^{0.74}\quad (5),\]
so that for $k(GV)\leq |V|$ it suffices to show that 
\[\frac{4}{5}|V|+3^{\frac{1}{2}(n-1)}|V|^{0.74}\leq |V|,\]
or, equivalently, 
\[5\cdot 3^{\frac{1}{2}(n-1)}\leq |V|^{0.26}=(25^n)^{0.26}=5^{0.52n}\quad (6).\]
For $n\geq 4$, this can indeed be checked to be true.\\
So finally let $n=3$. Any subgroup of $S_3$ contains at most two conjugacy classes of fixed
point free elements (namely, the two 3--cycles), and hence from \ref{lemb1} we see
that 
\begin{eqnarray*}
k(GV)&\leq&k(HV)+\mbox{ (number of conjugacy classes of elements of $G/N$ without fixed points on}\\
     &&\{V_1,V_2,V_3\})\cdot\max\{|C_{\scl(NV)}(g)|g\in G-\bigcup\limits_{x\in G}H^x\}\\
     &\leq&\frac{4}{5}|V|+2\cdot |V|^{0.74}\leq 15038<15625=|V|,
\end{eqnarray*}
and so the proof of the theorem is complete.
\end{beweis}

Note that techniques as in the above result also will work in may other interesting situations,
such as the ones that were left over by \cite[Theorem A]{riese-schmid}.\\

\section{Reducing the non--coprime $k(GV)$--problem}\label{sectionf}\\

In this section we present some ideas that might be helpful in dealing with the imprimitive case
of the following conjecture that has sometimes been called the non--coprime $k(GV)$--problem.\\

\begin{nconjecture}\label{conf1}
There is a universal constant $C$ such that the following holds:\\
Let $G$ be a finite $G$--module and $V$ be a finite faithful, completely reducible $G$--module. Then
\[k(GV)\leq C|V|\log_2|V|.\]
\end{nconjecture}
Note that our reduction results used in Section \ref{sectiond} are not always useful here, because
they require some knowledge (by induction) of $k(HV)$ where $H=N_G(V_1)$ for an imprimitivity
decomposition $V=V_1\oplus\ldots\oplus V_n$ of $V$, but in general there is no guarantee that $V$
is completely reducible as an $H$--module (although this is the case for small dimensions, see
\cite{guralnick}). Whenever $N\unlhd G$, however, then by Clifford $V$ is completely reducible
as an $N$--module, and so we can use \ref{lemb3} in combination with the 
following result of Guralnick and Magaard 
see \cite[Corollary 1]{guralnick-magaard}: If $G$ is a primitive permutation group on a set 
$\Omega$ of size $n$ and if $F^*(G)$ is not a product of alternating groups, then each nontrivial
element of $G$ fixes at least $\frac{4}{7}n$ elements of $\Omega$.\\
With this we can prove the following result.\\

\begin{nsatz}\label{theoremf2}
Let $f:\n\to\r$ be a function. Let $G$ be a finite group and $V$ be a finite $G$--module. Suppose that
$N\unlhd G$ and $V_N=V_1\oplus\ldots\oplus V_n$ for an $n\in\n$, where the $V_i$ are
$N$--modules. Assume further that $G/N$ primitively and faithfully permutes the $V_i$. Moreover
suppose that with $t_0=\max\{k(UV_1)\ |\ U\leq N/C_N(V_1)\}$ where
\[|N/C_N(V_1)|\leq\frac{\left(1-\frac{1}{|G/N|}\right)^\frac{14}{3n}f(|V|)^\frac{14}{3n}}{
2^\frac{14}{3}|V_1|t_0^\frac{8}{3}}\mbox{ and that }k(NV)\leq f(|V|).\]
Then one of the following holds:\\
(a) $k(GV)\leq f(|V|)$.\\
(b) $F^*(G/N)$ is a product of alternating groups (where $F^*(G/N)$ is the generalized Fitting
subgroup of $G/N$)
\end{nsatz}
\begin{beweis}
Clearly we may assume that $n>1$. Assume that $F^*(G/N)$ is not a product of alternating groups.
We have to show that (a) holds. By \cite[Corollary 1]{guralnick-magaard} we know that any
$g\in G-N$ fixes at most $\frac{4}{7}n$ of the $V_i$, and hence with \ref{lem1.2}(b) 
and \ref{lem2} we see that with $n_1:=|N/C_N(V_1)|$ we have
\[|C_{\scl(NV)}(g)|\leq t_0^{\frac{4}{7}n}\cdot (n_1|V_1|)^{\frac{3}{14}n}\]
for all $g\in G-N$.\\
Moreover, as $k(G/N)\leq 2^{n-1}$, with \ref{lemb3} we get
\begin{eqnarray*}
k(GV)&\leq&\frac{k(NV)}{|G/N|}+2(k(G/N)-1)\max\{|C_{\scl(NV)}(g)| |g\in G-N\}\\
     &\leq&\frac{f(|V|)}{|G/N|}+2^nt_0^{\frac{4}{7}n}|V_1|^{\frac{3}{14}n}n_1^{\frac{3}{14}n}\leq f(|V|)
\end{eqnarray*}
Thus by our hypothesis on $n_1$ we are done.
\end{beweis}

In view of \ref{conf1}, the following special case provides a reduction to primitive groups
in some situations.\\

\begin{ncor}\label{corf3}
Let $G$ be a finite group and $V$ be a finite $G$--module. Suppose that
$N\unlhd G$ and $V_N=V_1\oplus\ldots\oplus V_n$ for an $n\in\n$, where the $V_i$ are
$N$--modules. Assume further that $G/N$ primitively and faithfully permutes the $V_i$. 
Let $n\geq 5$ and $t_0$ be as in \ref{theoremf2}, and assume that 
\[|N/C_N(V_1)|\leq\frac{1}{50}C^\frac{14}{3n}\frac{|V_1|^\frac{11}{3}}{t_0^\frac{8}{3}}(\log_2|V|)^\frac{14}{3n}.\]
for some constant $C$. 
If $F^*(G/N)$ is not a product of alternating groups and $k(NV)\leq C|V|\log_2|V|$, then
\[k(GV)\leq C|V|\log_2|V|.\]
\end{ncor}
\begin{beweis}
Let $f(x)=Cx\log_2(x)$.
As $n\geq 5$, we see that $|G/N|\geq 5$ and $\frac{14}{3n}\leq\frac{14}{15}$ and hence
$\frac{1}{50}\leq\frac{\left(1-\frac{1}{5}\right)^\frac{14}{15}}{2^\frac{14}{3}}\leq
\frac{\left(1-\frac{1}{|G/N|}\right)^\frac{14}{3n}}{2^\frac{14}{3}}$, and thus our hypothesis
on $|N/C_N(V_1)|$ implies the one in \ref{theoremf2}. Hence by \ref{theoremf2} the assertion
follows.
\end{beweis}

Note that if $G$ is a minimal counterexample to \ref{conf1} and if $\char (V)\geq\dim V_1+2$,
then by the results in \cite{guralnick} we may assume that $t_0\leq C|V_1|\log_2|V_1|$, so that
by the hypothesis in \ref{corf3} may be replaced by the stronger condition that
\[|N/C_N(V_1)|\leq \frac{1}{50}C^{\frac{14}{3n}-\frac{8}{3}}|V_1|(\log_2|V|)^{\frac{14}{3n}-\frac{
8}{3}}.\]
This also yields Theorem D.\\

A similar result could be obtained with \ref{lemb2} in an obvious way, but we omit this here.\\
However, \ref{lemb1} gives a quite interesting result not involving $t_0$ (which in general can be
hard to control). \\

\begin{nsatz}\label{theorem7.4}
Let $G$ be a finite group and $V$ be a finite $G$--module. Suppose that $N\unlhd G$ and 
$V_N=V_1\oplus\ldots\oplus V_n$ for an $n\in\n$, where the $V_i$ are $N$--modules. Put
$H=N_G(V_1)$. Suppose that
\[k(HV)\leq C_1|V|\log_2|V|\]
for some constant $C_1$, and suppose that
\[|N/C_N(V_1)|\ \leq\  \frac{1}{4}\ (C_2-C_1)^\frac{2}{n}\ |V_1|\ (\log_2|V|)^\frac{2}{n}\]
for some constant $C_2$. Then
\[k(GV)\leq C_2|V|\log_2|V|.\]
\end{nsatz}
\begin{beweis}
Put $n_1=|N/C_N(V_1)|$. Then with \ref{lem1.2}(b) and \ref{lem2} for any $g\in G-\bigcup\limits_{
x\in G}H^x$ we have
\[|C_{\scl(NV)}(g)|\leq (n_1|V_1|)^\frac{n}{2},\]
and hence the assertion follows easily with \ref{lemb1}.
\end{beweis}

Note that the hypothesis on $H$ in an inductive proof of \ref{conf1} is satisfied whenever we know
that $H$ acts completely reducibly on $V$, which, for instance, by results of Guralnick
\cite{guralnick}, is the case whenever $\char (V)\geq\dim V+2$.

}


\begin{thebibliography}{99}
\bibitem{fulman-guralnick} {J. Fulman, R. Guralnick,} {\rm `Derangements in simple and primitive
groups'}, Groups, Combinatorics and geometry (Durham, 2001), World Scientific, New York 2003, 
99--121.
\bibitem{gallagher} {P. X. Gallagher,} {\rm `The number of conjugacy classes
of a finite group'}, {\it Math Z.} {\bf 118} (1970), 175--179.
\bibitem{gambini-weigel} {A. Gambini Weigel, T. S. Weigel,} {\rm `On the orders of primitive linear 
$p'$-groups'}, {\it Bull. Austral. Math. Soc.} {\bf 48} (1993), 495--521.
\bibitem{gluck-magaard} {D. Gluck, K. Magaard,} {\rm `Base sizes and regular orbits for coprime 
affine permutation groups'}, {\it J. London Math. Soc. (2)} {\bf 58} (1998), 603--618. 
\bibitem{gmrs} {D. Gluck, K. Magaard, U. Riese, P. Schmid,} {\rm `The solution of the
$k(GV)$--problem'}, {\it J. Algebra} {\bf 279} (2004), 694--719.
\bibitem{guralnick} {R. Guralnick,} {\rm `Small representations are completely reducible'}, 
{\it J. Algebra} {\bf 220} (1999), 531--541.
\bibitem{guralnick-magaard} {R. Guralnick, K. Magaard,} {\rm `On the minimal degree of a 
permutation group'}, {\it J. Algebra} {\bf 207} (1998), 127--145.
\bibitem{guralnick-tiep} {R. Guralnick, P. H. Tiep,} {\rm `The non--coprime $k(GV)$--problem'},
Preprint (2004).
\bibitem{guralnick-wan} {R. Guralnick, D. Wan,} {\rm `Bounds for fixed point free elements
in a transitive group and applications to curves over finite fields'}, {\it Israel J. Math.} 
{\bf 101} (1997), 255--287.
\bibitem{huppertcharacters} {B. Huppert,} {\rm `Character theory of finite groups'}, {deGruyter,}
Berlin, 1998.
\bibitem{isaacs} {M. I. Isaacs,} {\rm `Character theory of finite groups.'}, Academic Press, New 
York-London, 1976.
\bibitem{kellerkgv2} {T. M. Keller,} {\rm `The $k(GV)$--problem revisited'}, 
{\it J. Austral. Math. Soc.}, to appear.
\bibitem{kellerkgv1} {T. M. Keller,} {\rm `A new approach to the $k(GV)$--problem'}, 
{\it J. Austral. Math. Soc.} {\bf 75} (2003), 197--219.
\bibitem{kovacs-robinson} {L. G. Kov\'{a}cs, G. R. Robinson,} {\rm `On the number of conjugacy classes
of a finite group'}, {\it J. Algebra} {\bf 160} (1993), 441--460.
\bibitem{liebeck-pyber} {M. W. Liebeck, L. Pyber,} {\rm `Upper bounds for the number of conjugacy classes
of a finite group'}, {\it J. Algebra} {\bf 198} (1997), 538--562.
\bibitem{maroti} {A. Mar\'{o}ti,} {\rm `Bounding the number of conjugacy classes
of a permutation group'}, {\it J. Group Theory}, to appear. 
\bibitem{riese-schmid} {U. Riese, P. Schmid,} {\rm `Real vectors for linear groups and the
$k(GV)$--problem'}, {\it J. Algebra} {\bf 267} (2003), 725--755.
\bibitem{robinson-thompson} {G. R. Robinson, J. G. Thompson,} {\rm `On Brauer's $k(B)$--problem'},
{\it J. Algebra} {\bf 184} (1996), 1143--1160.
\end{thebibliography}
\end{document}